\documentclass[11pt, letterpaper, leqno]{amsart}
\usepackage{amsmath}
\usepackage{amsthm}
\usepackage{amsrefs}
\usepackage{lipsum}
\usepackage{amsopn}
\usepackage{amssymb}
\usepackage{enumerate} 
\usepackage{fullpage}
\numberwithin{equation}{section}
\usepackage{float}
\usepackage{graphicx}
\usepackage{titlesec}
\usepackage{fancyhdr}
\usepackage{fancyhdr}
\newtheorem{question}{Question}[section]

\newtheorem{corollary}{Corollary}[section]
\newtheorem{lemma}{Lemma}[section]
\newtheorem{theorem}{Theorem}[section]
\newtheorem{example}{Example}[section]

\newtheorem{definition}{Definition}[section]

\theoremstyle{definition}
\newtheorem{remark}{Remark}[section]

\DeclareMathOperator{\LOG}{Log}
\DeclareMathOperator{\Arg}{Arg}
\DeclareMathOperator{\adm}{adm}

\DeclareMathOperator{\dist}{dist}
\titleformat{\section}
  {\normalfont\scshape}{\thesection}{1em}{\centering}
\titleformat{\subsection}
{\normalfont\scshape}{\thesubsection}{1em}{\centering}

\begin{document}
\title{On the Hardy number of a domain in terms of harmonic measure and hyperbolic distance}

\author{Christina Karafyllia}  
\address{Department of Mathematics, Aristotle University of Thessaloniki, 54124, Thessaloniki, Greece}
\email{karafyllc@math.auth.gr}   
\thanks{I thank Professor D. Betsakos, my thesis advisor, for his advice during the preparation of this work, and the Onassis Foundation for the scholarship I receive during my Ph.D. studies.}

\fancyhf{}
\renewcommand{\headrulewidth}{0pt}

\fancyhead[RO,LE]{\small \thepage}
\fancyhead[CE]{\small On the Hardy number of a domain in terms of harmonic measure and hyperbolic distance}
\fancyhead[CO]{\small Christina Karafyllia} 
\fancyfoot[L,R,C]{}


\begin{abstract} Let $\psi $ be a conformal map on $\mathbb{D}$ with $ \psi \left( 0 \right)=0$ and let ${F_\alpha }=\left\{ {z \in \mathbb{D}:\left| {\psi \left( z \right)} \right| = \alpha } \right\}$ for $\alpha >0$. Denote by ${H^p}\left( \mathbb{D} \right)$ the classical Hardy space with exponent $p>0$ and by ${\tt h}\left( \psi \right)$ the Hardy number of $\psi$. Consider the limits 
\[L: = \mathop {\lim }\limits_{\alpha  \to  + \infty } \left( {{{\log {\omega _\mathbb{D}}{{\left( {0,{F_\alpha }} \right)}^{ - 1}}} \mathord{\left/
{\vphantom {{\log {\omega _D}{{\left( {0,{F_\alpha }} \right)}^{ - 1}}} {\log \alpha }}} \right.
\kern-\nulldelimiterspace} {\log \alpha }}} \right),\;\mu  := \mathop {\lim }\limits_{\alpha  \to  + \infty } \left( {{{d_\mathbb{D}}\left( {0,{F_\alpha }} \right)} \mathord{\left/
{\vphantom {{{d_D}\left( {0,{F_\alpha }} \right)} {\log \alpha }}} \right.
\kern-\nulldelimiterspace} {\log \alpha }} \right),\]
where $\omega _\mathbb{D}\left( {0,{F_\alpha }} \right)$ denotes the harmonic measure at $0$ of $F_\alpha $ and $d_\mathbb{D} {\left( {0,{F_\alpha }} \right)}$ denotes the hyperbolic distance between $0$ and $F_\alpha$ in $\mathbb{D}$. We study a problem posed by P. Poggi-Corradini. What is the relation between 
$L$, $\mu$  and ${\tt h}\left( \psi \right)$? Motivated by the result of Kim and Sugawa that ${\tt h}\left( \psi \right)= \mathop {\liminf }\limits_{\alpha  \to  + \infty } \left( {{{\log {\omega _\mathbb{D}}{{\left( {0,{F_\alpha }} \right)}^{ - 1}}} \mathord{\left/
			{\vphantom {{\log {\omega _D}{{\left( {0,{F_\alpha }} \right)}^{ - 1}}} {\log \alpha }}} \right.
			\kern-\nulldelimiterspace} {\log \alpha }}} \right)$, we show that ${\tt h}\left( \psi \right) = \mathop {\liminf }\limits_{\alpha  \to  + \infty } \left( {{{d_\mathbb{D}}\left( {0,{F_\alpha }} \right)} \mathord{\left/
		{\vphantom {{{d_D}\left( {0,{F_\alpha }} \right)} {\log \alpha }}} \right.
		\kern-\nulldelimiterspace} {\log \alpha }} \right)$. We also provide conditions for the existence of $L$ and $\mu$ and for the equalities $L=\mu={\tt h}\left( \psi \right)$.
Poggi-Corradini proved that $\psi  \notin {H^{\mu}}\left( \mathbb{D} \right)$ for a wide class of conformal maps $\psi$. We present an example of $\psi$ such that $\psi  \in {H^\mu {\left( \mathbb{D} \right)} }$.
\end{abstract}

\maketitle
\section{Introduction}\label{section1}

We study the Hardy number of a domain in terms of harmonic measure and hyperbolic distance. For a domain $D$, a point $z \in D$ and a Borel subset $E$ of $\overline D $, let ${\omega _D}\left( {z,E} \right)$ denote the harmonic measure at $z$ of $E$ with respect to the component of $D \backslash {E}$ containing $z$. The function ${\omega _D}\left( { \cdot ,E} \right)$ is exactly the solution of the generalized Dirichlet problem with boundary data $\varphi  = {1_E}$   (see \cite[ch. 3]{Ahl}, \cite[ch. 1]{Gar} and \cite[ch. 4]{Ra}). The hyperbolic distance between two points $z,w$ in the unit disk $\mathbb{D}$ (see \cite[ch. 1]{Ahl}, \cite[p. 11-28]{Bea}) is defined by 
\[{d_\mathbb{D}}\left( {z,w} \right) = \log \frac{{1 + \left| {\frac{{z - w}}{{1 - z\bar w}}} \right|}}{{1 - \left| {\frac{{z - w}}{{1 - z\bar w}}} \right|}}.\]
The hyperbolic distance can be defined on any simply connected domain $D \ne \mathbb{C}$ as follows: If $f$ is a Riemann map of $\mathbb{D}$ onto $D$ and $z,w \in D$, then 
${d_D}\left( {z,w} \right) = {d_\mathbb{D}}\left( {{f^{ - 1}}\left( z \right),{f^{ - 1}}\left( w \right)} \right)$.
Also, for a set $E \subset D$, we define ${d_D}\left( {z,E} \right): = \inf \left\{ {{d_D}\left( {z,w} \right):w \in E} \right\}$.

The Hardy space with exponent $p$, $p>0$, and norm ${\left\|  \cdot  \right\|_p}$ (see \cite[p. 1-2]{Du}, \cite[p. 435-441]{Gar}) is defined to be
\[{H^p}\left( \mathbb{D} \right) = \left\{ {f \in H\left( \mathbb{D} \right):\left\| f \right\|_p^p = \mathop {\sup }\limits_{0 < r < 1} \int_0^{2\pi } {{{\left| {f\left( {r{e^{i\theta }}} \right)} \right|}^p}d\theta  <  + \infty } } \right\},\]
where $H\left( \mathbb{D} \right)$ denotes the family of all holomorphic functions on $\mathbb{D}$.  The fact that a function $f$ belongs to ${H^p}\left( \mathbb{D} \right)$ imposes a restriction on the growth of $f$ and this restriction is stronger as $p$ increases. If $\psi$ is a conformal map on $\mathbb{D}$, then $\psi \in {H^p}\left( \mathbb{D} \right)$ for all $p<1/2$ (\cite[p. 50]{Du}). 

\begin{figure}[H]
	\begin{center}
		\includegraphics[scale=0.5]{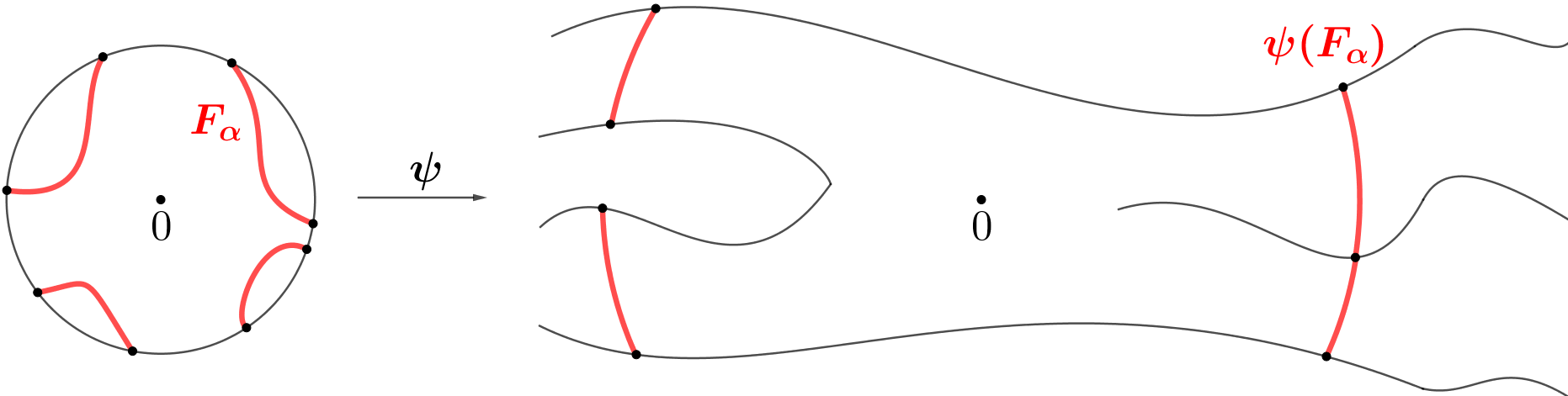}
		\caption{The conformal map $\psi$ on $\mathbb{D}$ and the sets $F_\alpha, \psi \left( F_\alpha \right)$.}
		\label{ra}
	\end{center}
\end{figure}

Hereinafter, $\psi $ is a conformal map on $\mathbb{D}$ with $ \psi \left( 0 \right)=0$ and ${F_\alpha } = \left\{ {z \in \mathbb{D}:\left| {\psi \left( z \right)} \right| = \alpha } \right\}$ for $\alpha >0$ (see Fig. \ref{ra}). The number ${\tt h}\left( \psi \right) \in \left[ {{1 \mathord{\left/
 {\vphantom {1 2}} \right.
 \kern-\nulldelimiterspace} 2}, + \infty } \right]$ which is given by
\[{\tt h}\left( \psi \right) = \sup \left\{ {p > 0:\psi  \in {H^p}\left( \mathbb{D} \right)} \right\},\]
is called the Hardy number of $\psi$ and  was first introduced by Hansen in \cite{Ha}. Note that if $D$ is a simply connected domain, we say $D \in {H^p}\left( \mathbb{D} \right)$ if there is a Riemann map $\psi$ of $\mathbb{D}$ onto $D$ such that $\psi \in {H^p}\left( \mathbb{D} \right)$. Any other Riemann map onto $D$ is also in ${H^p}\left( \mathbb{D} \right)$, and hence the Hardy number of $D$ is well-defined by setting ${\tt h}\left( D \right)={\tt h}\left( \psi \right)$. A classical problem in geometric function theory is to find the Hardy number of a domain by looking at its geometric properties (see e.g. \cite{Ba}, \cite{Rat}). Hansen studied the number by using Ahlfors' distortion theorem and he described it in terms of geometric quantities for starlike and spiral-like domains \cite{Han}. In \cite{Es} Ess{\'e}n gave a description of ${\tt h}\left( \psi \right)$ in terms of harmonic measures and obtained almost necessary and sufficient conditions for ${\tt h}\left( \psi \right)$ in terms of capacity. Poggi-Corradini \cite{Co1} studied the range domains $D$ of univalent K{\oe}nigs functions (see also \cite{Co2}) and found that the number ${\tt h}\left( D \right)$ can be described in terms of the essential norm of the associated composition operators. Finally, based on Ess{\'e}n' s main lemma \cite{Es}, Kim and Sugawa \cite{Kim} proved that 

\begin{equation}\label{sug}
{\tt h}\left( \psi \right)=\mathop {\liminf}\limits_{\alpha  \to  + \infty } \frac{{\log {\omega _{\psi \left( \mathbb{D} \right)}}{{\left( {0,{\psi \left( F_\alpha \right)}} \right)}^{ - 1}}}}{{\log \alpha }} = \mathop {\liminf}\limits_{\alpha  \to  + \infty } \frac{{\log {\omega _\mathbb{D}}{{\left( {0,{F_\alpha }} \right)}^{ - 1}}}}{{\log \alpha }}.
\end{equation}
In Section \ref{section4} we express ${\tt h}\left( \psi \right)$ in terms of hyperbolic distance by proving the following theorem.

\begin{theorem}\label{hn} Let $\psi $ be a conformal map on $\mathbb{D}$ with $ \psi \left( 0 \right)=0$ and let ${F_\alpha } = \left\{ {z \in \mathbb{D}:\left| {\psi \left( z \right)} \right| = \alpha } \right\}$ for $\alpha >0$. If ${\tt h}\left( \psi \right)$ denotes the Hardy number of $\psi$, then
\[{\tt h}\left( \psi \right) = \mathop {\liminf}\limits_{\alpha  \to  + \infty } \frac{{{d_\mathbb{D}}\left( {0,{F_\alpha }} \right)}}{{\log \alpha }}.\]
\end{theorem}

\begin{figure}[H]
	\begin{center}
		\includegraphics[scale=0.5]{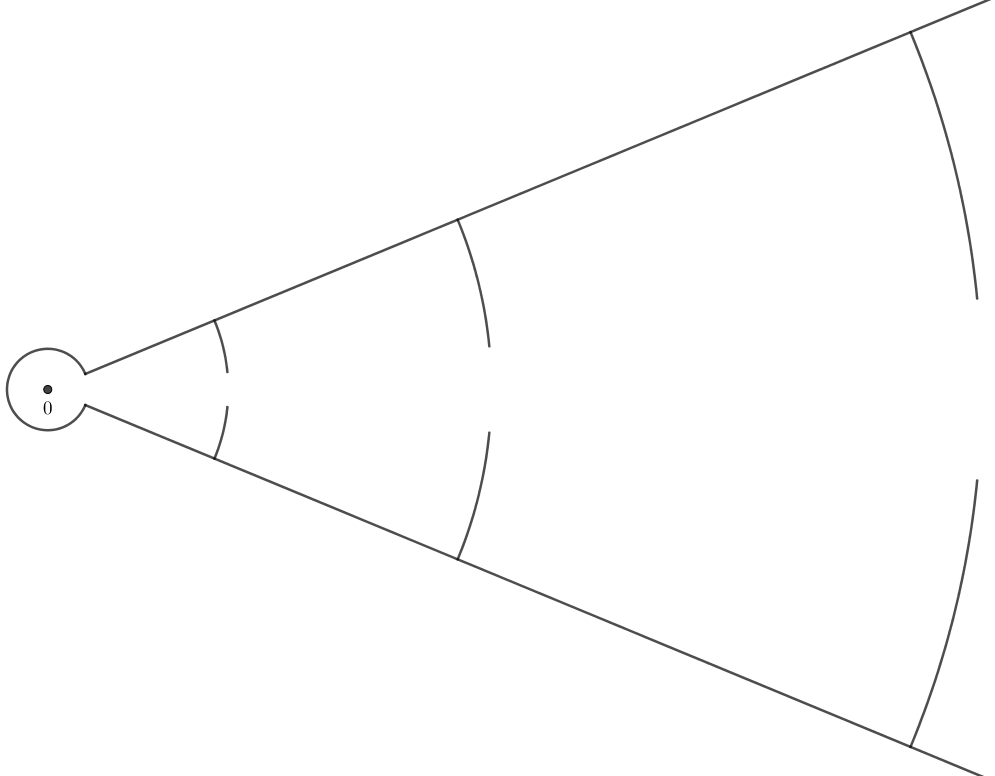}
		\caption{}
		\label{topos1}
	\end{center}
\end{figure}

Harmonic measure and hyperbolic distance are both conformally invariant and several Euclidean estimates are known about them. Thus, expressing the ${H^p}\left( \mathbb{D} \right)$-norms of a conformal map $\psi$ on $\mathbb{D}$ in terms of harmonic measure and hyperbolic distance, we are able to obtain information about the growth of the function by looking at the geometry of its image region $\psi \left( {\mathbb{D}} \right)$. In \cite[p. 10]{Co} Poggi-Corradini proved that the Beurling-Nevanlinna projection theorem \cite[p. 43-44]{Ahl} implies that for every $\alpha>0$,
\[{\omega _\mathbb{D}}\left( {0,{F_\alpha }} \right) \ge \frac{2}{\pi }{e^{ - {d_\mathbb{D}}\left( {0,{F_\alpha }} \right)}}\]
and he stated the question \cite[p. 36]{Co} whether the opposite inequality is also true for some positive constant. In \cite{Ka} we proved that the answer is negative and only under additional assumptions involving the  geometry of the domain $\psi \left( \mathbb{D} \right)$ it can be positive. However, the situation changes when we study integrals of the quantities stated above. In \cite[p. 33]{Co} and \cite[p. 502-503]{Co2} Poggi-Corradini proved that
\begin{equation}\label{isod}
\psi  \in {H^p}\left( \mathbb{D} \right) \Leftrightarrow \int_0^{ + \infty } {{\alpha ^{p - 1}}{\omega _{\mathbb{D}}}\left( {0,{F_\alpha }} \right)d\alpha }  <  + \infty.
\end{equation}
Answering a question he stated in \cite[p. 36]{Co}, we proved in \cite{Kar} that
\begin{equation}\label{isod1}
\psi  \in {H^p}\left( {\mathbb{D}} \right) \Leftrightarrow \int_0^{ + \infty } {{\alpha ^{p - 1}}{e^{ - {d_{\mathbb{D}}}\left( {0,{F_\alpha }} \right)}}d\alpha }  <  + \infty.
\end{equation}
If we rewrite the integrands of conditions (\ref{isod}) and (\ref{isod1}), we take respectively,
\[{\alpha ^{p - 1}}{\omega _\mathbb{D}}\left( {0,{F_\alpha }} \right) = {\alpha ^{p - 1 - {{\log {\omega _\mathbb{D}}{{\left( {0,{F_\alpha }} \right)}^{ - 1}}} \mathord{\left/
 {\vphantom {{\log {\omega _\mathbb{D}}{{\left( {0,{F_\alpha }} \right)}^{ - 1}}} {\log \alpha }}} \right.
 \kern-\nulldelimiterspace} {\log \alpha }}}}\]
and
\[{\alpha ^{p - 1}}{e^{ - {d_{\mathbb{D}}}\left( {0,{F_\alpha }} \right)}}={\alpha ^{p - 1 - {{{d_\mathbb{D}}\left( {0,{F_\alpha }} \right)} \mathord{\left/
				{\vphantom {{{d_{\mathbb{D}}}\left( {0,{F_\alpha }} \right)} {\log \alpha }}} \right.
				\kern-\nulldelimiterspace} {\log \alpha }}}}.\]
Poggi-Corradini noticed that if the limit $L: = \mathop {\lim }\limits_{\alpha  \to  + \infty } \left( {{{\log {\omega _\mathbb{D}}{{\left( {0,{F_\alpha }} \right)}^{ - 1}}} \mathord{\left/
 {\vphantom {{\log {\omega _D}{{\left( {0,{F_\alpha }} \right)}^{ - 1}}} {\log \alpha }}} \right.
 \kern-\nulldelimiterspace} {\log \alpha }}} \right)$ exists then the ratio ${{{\log {\omega _\mathbb{D}}{{\left( {0,{F_\alpha }} \right)}^{ - 1}}} \mathord{\left/
  {\vphantom {{\log {\omega _D}{{\left( {0,{F_\alpha }} \right)}^{ - 1}}} {\log \alpha }}} \right.
  \kern-\nulldelimiterspace} {\log \alpha }}}$ determines the Hardy number of $\psi$. In fact, by (\ref{isod}) we deduce that if $p<L$ then $\psi  \in {H^p}\left( \mathbb{D} \right)$ and if $p>L$, $\psi  \notin {H^p}\left( \mathbb{D} \right)$. Similarly, if the limit $\mu  := \mathop {\lim }\limits_{\alpha  \to  + \infty } \left( {{{d_\mathbb{D}}\left( {0,{F_\alpha }} \right)} \mathord{\left/
{\vphantom {{{d_D}\left( {0,{F_\alpha }} \right)} {\log \alpha }}} \right.
\kern-\nulldelimiterspace} {\log \alpha }} \right)$ exists then by (\ref{isod1}) we infer that if $p<\mu$ then $\psi  \in {H^p}\left( \mathbb{D} \right)$ and if $p>\mu$ then $\psi  \notin {H^p}\left( \mathbb{D} \right)$. So, the ratio ${{{d_\mathbb{D}}\left( {0,{F_\alpha }} \right)} \mathord{\left/
{\vphantom {{{d_D}\left( {0,{F_\alpha }} \right)} {\log \alpha }}} \right.
\kern-\nulldelimiterspace} {\log \alpha }}$ determines the Hardy number of $\psi$. However, it is not clear whether  $\psi  \in {H^p}\left( \mathbb{D} \right)$ when $\mu$ (or $L$) is finite and $p=\mu$ (or $p=L$). Poggi-Corradini proved (see \cite[p. 37-38]{Co} and \cite[p. 503-504]{Co2}) that $\psi  \notin {H^{\mu}}\left( \mathbb{D} \right)$ for a wide class of conformal maps $\psi$ which he calls ``sector-like". But, could this result be generalized for every simply connected domain? In Section \ref{section5}, we answer this question by constructing the simply connected domain of Fig. \ref{topos1} so that, if $\psi$ is the corresponding Riemann map, then $\psi  \in {H^\mu {\left( \mathbb{D} \right)} }$. The reasons, which led us to construct this particular domain, are stated at the beginning of Section \ref{section5}.

\begin{example}\label{ex} There exists a conformal map $\psi$ on $\mathbb{D}$ such that $\mu$ exists and $\psi  \in {H^\mu {\left( \mathbb{D} \right)} }$.
\end{example} 

Therefore, when $\mu$ (or $L$) is finite, the case $p=\mu$ (or $p=L$) depends on the way the ratio approaches the limit $\mu$ (or $L$). Finally, to complete the study of these limits, it is reasonable to examine the connection between $\mu$ and $L$. So, in Section \ref{section4}, we prove the following results. 

\begin{theorem}\label{thL} Suppose that $\mu $ exists. Then $L$ exists and $L= \mu$. 
\end{theorem}

\begin{corollary}\label{coth} $\mu  =  + \infty$ if and only if $L =  + \infty$. 
\end{corollary}

Let  $N\left( \alpha \right) \in \mathbb{N} \cup \left\{ { + \infty } \right\}$ denote the number of components of $F_\alpha$ for $\alpha>0$ and $F_\alpha ^i$ denote each of these components for $i=1,2, \ldots, N\left( \alpha  \right)$. Since $\max \left\{ {{\omega _\mathbb{D}}\left( {0,F_\alpha ^i} \right):i \in \left\{ {1,2, \ldots N\left( \alpha  \right)} \right\}} \right\}$ exists, as we prove in Section \ref{section3}, we denote by $F_\alpha ^ *$ a component of $F_\alpha$ such that
\[{\omega _\mathbb{D}}\left( {0,F_\alpha ^ * } \right) = \max \left\{ {{\omega _\mathbb{D}}\left( {0,F_\alpha ^i} \right):i \in \left\{ {1,2, \ldots, N\left( \alpha  \right)} \right\}} \right\}.\]

\begin{theorem}\label{cond} Suppose that $L$ exists. Then $\mu $ exists if and only if 
\begin{equation}\label{condit}
\mathop {\lim \sup }\limits_{\alpha  \to  + \infty } \frac{{\log {\omega _\mathbb{D}}{{\left( {0,F_\alpha ^ * } \right)}^{ - 1}}}}{{\log \alpha }} = L.
\end{equation}
In case $\mu $ exists then $\mu=L $.
\end{theorem}

\begin{corollary}\label{cor} If $L$ exists and $\mathop {\lim }\limits_{\alpha  \to  + \infty } \frac{\log{N\left( \alpha  \right)}}{{\log \alpha }} = 0$ then $\mu $ exists and $\mu=L $.
\end{corollary}

Note that the condition of the corollary above is more geometric and easy to check but it is not clear if it is necessary and sufficient. On the other hand, the condition (\ref{condit}) of Theorem \ref{cond} is necessary and sufficient but not so easy to handle. So, we state the following question. 

\begin{question} Can we replace the condition {\rm (\ref{condit})} by a more geometric condition or, maybe, is the condition {\rm (\ref{condit})} true for every simply connected domain?
\end{question}

In Section \ref{section2} we introduce some preliminary  results and notions such as the domain decomposition method studied by N. Papamichael and N.S. Stylianopoulos \cite{Pab}, the extremal length and its connection with the harmonic measure. In Section \ref{section3} we present some lemmas required for the proofs of Section \ref{section4}. In Section \ref{section4} we prove Theorems \ref {hn}, \ref{thL} and \ref{cond} and Corollaries \ref{coth} and \ref{cor}. Finally, in Section \ref{section5} we present the conformal map of the Example \ref{ex}.

\section{Preliminaries}\label{section2}

We first state a theorem proved by  Poggi-Corradini in \cite[p. 37]{Co} and \cite[p.134]{Co1}.   
\begin{theorem}\label{Corra}
Let $\psi$ be a conformal map on $\mathbb{D}$ and, for $\alpha>0$, let ${F_\alpha } = \left\{ {z \in \mathbb{D}:\left| {\psi \left( z \right)} \right| = \alpha } \right\}$. 
\begin{enumerate}[\rm (i)]
\item If $S = \mathop {\lim \sup }\limits_{\alpha  \to  + \infty } \frac{{{d_\mathbb{D}}\left( {0,{F_\alpha }} \right)}}{{\log \alpha }}<+\infty $, then:
	\begin{enumerate}
		\item $S < p <  + \infty  \Rightarrow \psi  \notin {H^p}\left( \mathbb{D} \right)$
		\item ${\alpha ^{S - 1 - {{{d_\mathbb{D}}\left( {0,{F_\alpha }} \right)} \mathord{\left/
						{\vphantom {{{d_{\mathbb{D}}}\left( {0,{F_\alpha }} \right)} {\log \alpha }}} \right.
						\kern-\nulldelimiterspace} {\log \alpha }}}}$ not integrable at infinity $\Rightarrow \psi  \notin {H^S}\left( \mathbb{D} \right)$. 
	\end{enumerate}
\item If $I = \mathop {\liminf}\limits_{\alpha  \to  + \infty } \frac{{{d_\mathbb{D}}\left( {0,{F_\alpha }} \right)}}{{\log \alpha }}$, then $I \ge {1 \mathord{\left/
 {\vphantom {1 2}} \right.
 \kern-\nulldelimiterspace} 2}$ and
\[0 < p <  I  \Rightarrow \psi  \in {H^p}\left( \mathbb{D} \right).\]
\end{enumerate}
In particular, if $S=I=\mu$ then $\mu  = \tt h\left( {\psi } \right)$.
\end{theorem}

\subsection{Extremal length}

Another  conformally invariant quantity, which is related to the harmonic measure, is the extremal length. We present the definition and the properties we need as they are stated in \cite[ch. 4]{Ahl}, \cite[p. 361-385]{Be}, \cite[ch. 7]{Fu}, \cite[ch. 4]{Gar} and \cite[ch. 2]{Oh}.

\begin{definition} Let $\left\{ C \right\}$ be a family of curves and $\rho \left( z \right) \ge 0$ be a measurable function defined in $\mathbb{C}$. We say $\rho \left( z \right)$ is admissible for $\left\{ C \right\}$ and denote by $\rho  \in \adm \left\{ C \right\}$, if for every rectifiable $C \in \left\{ C \right\}$, the integral $\int_C {\rho \left( z \right)\left| {dz} \right|} $ exists and $1 \le \int_C {\rho \left( z \right)\left| {dz} \right|}  \le  + \infty $. The extremal length of $\left\{ C \right\}$, $\lambda  \left\{ C \right\}$, is defined by
	\[\frac{1}{{\lambda \left\{ C \right\}}} = \mathop {\inf }\limits_{\rho  \in {\rm adm}\left\{ C \right\}} \int \int {{\rho ^2}\left( z \right)dxdy}. \]
\end{definition}

Note that if all curves of $\left\{ C \right\}$ lie in a domain $D$, we may take ${\rho \left( z \right) = 0}$ outside $D$. The conformal invariance is an immediate consequence of the definition (see \cite[p. 90]{Fu}). As a typical example (see \cite[p. 366]{Be}, \cite[p. 131]{Gar}), we mention the case in which $R$ is a rectangle with sides of length $a$ and $b$ and $\left\{ C \right\}$ is the family of curves in $R$ joining the opposite sides of length $a$. Then $\lambda \left\{ C \right\} = \frac{b}{a}$. Next we state two basic properties of extremal length that we will need (see \cite[p. 54-55]{Ahl}, \cite[p. 363]{Be}, \cite[p. 91]{Fu}, \cite[p. 134-135]{Gar}, \cite[p. 79]{Oh}).

\begin{theorem}\label{el} If $\left\{ {C'} \right\} \subset \left\{ C \right\}$ or every $C' \in \left\{ {C'} \right\}$ contains a $C \in \left\{ {C} \right\}$, then $\lambda \left\{ C \right\} \le \lambda \left\{ {C'} \right\}.$
\end{theorem}

\begin{theorem}[The serial rule]\label{sr} Let $\left\{ {B_n} \right\}$ be mutually disjoint Borel sets and each $C_n \in \left\{ {C_n} \right\}$ be in $B_n$. If $ \left\{ {C} \right\}$ is a family of curves such that each $C$ contains at least one $C_n$ for every $n$, then
	\[\lambda \left\{ C \right\} \ge \sum\limits_n {\lambda \left\{ {{C_n}} \right\}}. \] 
\end{theorem} 

Sometimes it is more convenient to use the more special notion of extremal distance. Let $D$ be a plane domain and $E_1,E_2$ be two disjoint closed sets on $\partial D$. If $\left\{ {C} \right\}$ is the family of curves in $D$ joining $E_1$ to $E_2$, then the extremal length $\lambda_D \left\{ C \right\}$ is called the extremal distance between $E_1$ and $E_2$ with respect to $D$ and is denoted by ${\lambda _D}\left( {{E_1},{E_2}} \right)$.

\subsection{Domain decomposition method}

In case of quadrilaterals, the opposite inequality in the serial rule has been studied by Papamichael and Stylianopoulos by means of a domain decomposition method for approximating the conformal modules of long quadrilaterals (see \cite{Pab}). Before stating the theorems we need, we present the required notation. 

\begin{figure}[H]
	\begin{center}
		\includegraphics[scale=0.5]{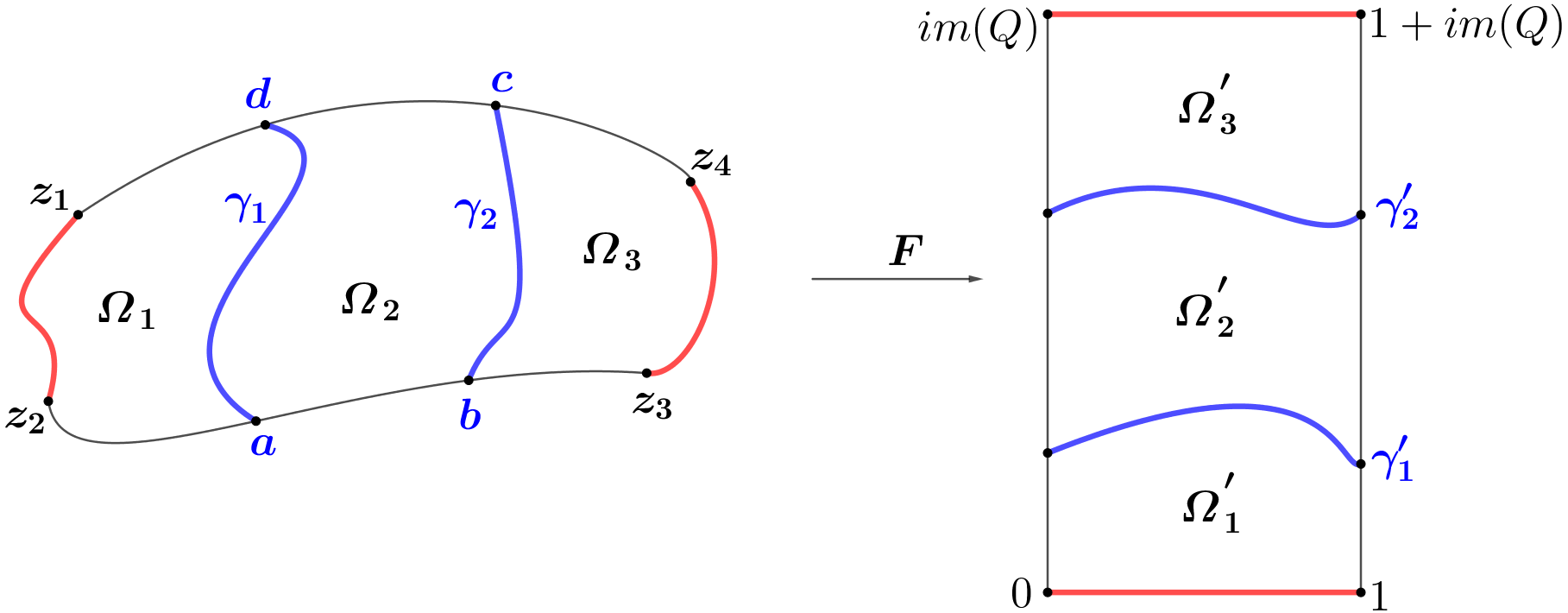}
		\caption{The subdivision of $\Omega$ into $\Omega_1,\Omega_2,\Omega_3$ and the conformal map $F:Q \to Q'$.}
		\label{dd}
	\end{center}
\end{figure}

Let $\Omega $ be a Jordan domain in $\mathbb{C}$ and consinder a system consisting of $\Omega $ and four distinct points $z_1,z_2,z_3,z_4$ in counterclockwise order on its boundary $\partial \Omega$. Such a system is said to be a quadrilateral $Q$ and is denoted by
\[Q: = \left\{ {\Omega ;{z_1},{z_2},{z_3},{z_4}} \right\}.\] 
The conformal module $m\left( Q \right)$ of $Q$ is the unique number for which $Q$ is conformally equivalent to the rectangular quadrilateral
\[Q': = \left\{ {R_{m\left( Q \right)} ;0,1,1+m\left( Q \right)i,m\left( Q \right)i} \right\},\] 
where $R_{m\left( Q \right)}=\left\{ {x + yi:0 < x < 1,0 < y < m\left( Q \right)} \right\}$ (see Fig. \ref{dd}). Note that $m\left( Q \right)$ is conformally invariant and it is equal to the extremal distance between the boundary arcs $\left( {z_1},{z_2} \right)$ and $\left( {z_3},{z_4} \right)$ of $\Omega$. So, $\Omega$ and $Q: = \left\{ {\Omega ;{z_1},{z_2},{z_3},{z_4}} \right\}$ will denote respectively the original domain and the corresponding quadrilateral. Moreover, ${\Omega _1},{\Omega _2}, \ldots ,$ and ${Q _1},{Q _2}, \ldots ,$ will denote the principle subdomains and corresponding component quadrilaterals of the decomposition under considerartion. Now consider the situation of Fig. \ref{dd}, where the decomposition of $Q: = \left\{ {\Omega ;{z_1},{z_2},{z_3},{z_4}} \right\}$ is defined by two non-intersecting arcs $\gamma_1,\gamma_2$ that join respectively two distinct points $a$ and $b$ on the boundary arc $\left( {{z_2},{z_3}} \right)$ to two points $d$ and $c$ on the boundary arc $\left( {{z_4},{z_1}} \right)$. These two arcs subdivide $\Omega$ into three non-intersecting subdomains denoted by ${\Omega _1},{\Omega _2}$ and ${\Omega _3}$. In addition, the arc $\gamma_1$ subdivides $\Omega$ into $\Omega_1$ and another subdomain denoted by ${\Omega _{2,3}}$, i.e. we take
\[{\overline \Omega  _{2,3}} = {\overline \Omega  _2} \cup {\overline \Omega  _3}.\]
Similarly, we say that $\gamma_2$ subdivides $\Omega$ into $\Omega_{1,2}$ and $\Omega  _3$, i.e. we take
\[{\overline \Omega  _{1,2}} = {\overline \Omega  _1} \cup {\overline \Omega  _2}.\]
Finally, we use the notations $Q_1,Q_2,Q_3,Q_{1,2}$ and $Q_{2,3}$ to denote, respectively, the quadrilaterals corresponding to the subdomains $\Omega_1,\Omega_2,\Omega_3,\Omega_{1,2}$ and $\Omega_{2,3}$, i.e.
\[Q_1: = \left\{ {\Omega_1 ;{z_1},{z_2},a,d} \right\},\;Q_2: = \left\{ {\Omega_2 ;d,a,b,c} \right\},\;Q_3: = \left\{ {\Omega_3 ;c,b,z_3,z_4} \right\}\]
and
\[Q_{1,2}: = \left\{ {\Omega_{1,2} ;{z_1},{z_2},b,c} \right\},\;Q_{2,3}: = \left\{ {\Omega_{2,3} ;d,a,z_3,z_4} \right\}.\]
The following theorems were proved by Papamichael and Stylianopoulos in \cite[p. 142-145]{Pab}.

\begin{theorem}\label{ddm1} Consider the decomposition and the notations illustrated in Fig. \ref{dd}. With the terminology  defined above, we have
	\[\left| {m\left( Q \right) - \left( {m\left( {{Q_{1,2}}} \right) + m\left( {{Q_{2,3}}} \right) - m\left( {{Q_{2}}} \right)} \right)} \right| \le 2.71{e^{ - \pi m\left( {{Q_2}} \right)}},\]
	provided that $m\left( {{Q_2}} \right) \ge 3$.
\end{theorem}

\begin{figure}[H]
	\begin{center}
		\includegraphics[scale=0.7]{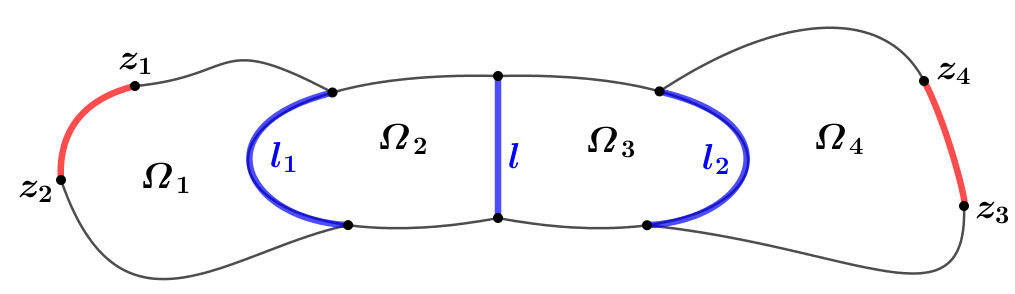}
		\caption{The decomposition of Theorem \ref{ddm}.}
		\label{dd1}
	\end{center}
\end{figure}

\begin{theorem}\label{ddm} Consider a quadrilateral $Q: = \left\{ {\Omega ;{z_1},{z_2},{z_3},{z_4}} \right\}$ of the form illustrated in Fig. \ref{dd1} and assume that the defining domain $\Omega  $ can be decomposed by means of a straight line crosscut $l$ and two other crosscuts $l_1$ and $l_2$ into four subdomains $\Omega_1$, $\Omega_2$, $\Omega_3$ and $\Omega_4$, so that $\Omega_3$ is the reflection in $l$ of $\Omega_2$. Then, for the decomposition of $Q$ defined by $l$,
\[0 \le m\left( Q \right) - \left( {m\left( {{Q_{1,2}}} \right) + m\left( {{Q_{3,4}}} \right)} \right) \le 5.26{e^{ - 2\pi m\left( {{Q_2}} \right)}},\]
provided that $m\left( {{Q_2}} \right) \ge 1.5$.
\end{theorem}

\begin{remark}\label{re} Papamichael and Stylianopoulos proved Theorems \ref{ddm1} and \ref{ddm} in case $\Omega $ is a Jordan domain. However, it follows from the proof that they are still valid if $\Omega $ is a simply connected domain and its boundary sets $\left( {{z_1},{z_2}} \right)$ and $\left( {{z_3},{z_4}} \right)$ are arcs of prime ends.
\end{remark}

\subsection{Harmonic measure}

Next we state a version of the Beurling-Nevanlinna projection theorem (see \cite[p. 43-44]{Ahl}, \cite[p. 43]{Be}, \cite[p. 105]{Gar} and \cite[p. 120]{Ra}) which gives us a relation between the harmonic measure of a closed and connected set in $\mathbb{D}$ and the harmonic measure of its circular projection on the negative radius.

\begin{theorem}[Beurling-Nevanlinna projection theorem]\label{proj} Let $E \subset {{\overline {\mathbb{D}}\backslash \left\{ 0 \right\}}}$ be a closed connected set intersecting the unit circle. Let ${E^ * } = \left\{ { - \left| z \right|:z \in E} \right\} = \left( { - 1,} \right.\left. { - {r_0}} \right]$, where ${r_0} = \min \left\{ {\left| z \right|:z \in E} \right\}$. Then, for $0 \le x < 1$,
\[{\omega _{\mathbb{D}}}\left( {x,E} \right) \ge {\omega _{\mathbb{D} }}\left( {x,{E^ * }} \right) = \frac{2}{\pi }\arcsin \frac{{\left( {1 - {r_0}} \right) \left( {1 - x} \right)}}{{\left( {1 + {r_0}} \right) \left( {1 + x} \right)}}.\]
\end{theorem}

Harmonic measure increases as the domain, in which it is defined, extends (see \cite[p. 102]{Ra}).

\begin{theorem}\label{pro} Let $D_1$, $D_2$ be simply connected domains such that ${D_1} \subset {D_2}$ and $B$ be a Borel subset of $\partial {D_1} \cap \partial {D_2}$. Then, for $z \in D_1$,
\[{\omega _{{D_1}}}\left( {z,B} \right) \le {\omega _{{D_2}}}\left( {z,B} \right).\]
\end{theorem}

Let $D$ be a bounded simply connected domain, $E$ be an arc on $\partial D$ and $z_0 \in D$. Consider all Jordan arcs $\sigma  \subset D$ that join $z_0$ to $\partial D\backslash E$ and define
\[{\lambda _D}\left( {{z_0},E} \right) = \mathop {\sup }\limits_\sigma  {\lambda _{D\backslash \sigma }}\left( {\sigma ,E} \right),\]
where the supremum is taken over all such Jordan arcs. Then the following theorem gives a relation between ${\omega _D}\left( {{z_0},E} \right)$ and ${\lambda _D}\left( {{z_0},E} \right)$ (see \cite[p. 368-371]{Be}, \cite[p. 144-146]{Gar}). 

\begin{theorem}\label{mar} Let $D$ be a bounded simply connected domain, $E$ be an arc on $\partial D$ and $z_0 \in D$. Then
	\[{e^{ - \pi {\lambda _D}\left( {{z_0},E} \right)}} \le {\omega _D}\left( {{z_0},E} \right) \le \frac{8}{\pi }{e^{ - \pi {\lambda _D}\left( {{z_0},E} \right)}}.\]
\end{theorem}

\section{Auxilary lemmas}\label{section3}

\begin{lemma}\label{geodesic}
	Let $\Gamma $ be the hyperbolic geodesic joining two points $z_1,z_2 \in \partial{\mathbb{D}}$ in $\mathbb{D}$. Then
	\[{e^{ - {d_{\mathbb{D}}}\left( {0,\Gamma } \right)}} \le {\omega _{\mathbb{D}}}\left( {0,\Gamma } \right) \le \frac{4}{\pi }{e^{ - {d_{\mathbb{D}}}\left( {0,\Gamma } \right)}}.\]
\end{lemma}

\begin{figure}[H] 
	\begin{center}
		\includegraphics[scale=0.7]{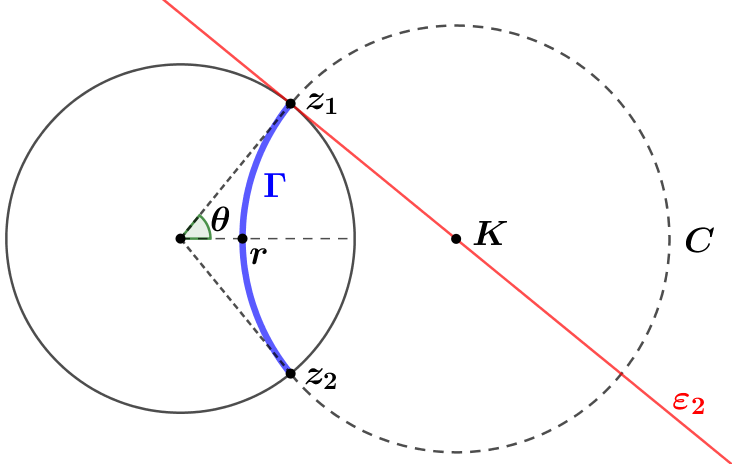}
		\caption{}
		\label{geodes}
	\end{center}
\end{figure}

\proof Without loss of generality, let ${z_1} = {e^{i\theta }}$, ${z_2} = {e^{-i\theta }}$ for some $\theta  \in \left( {0,\frac{\pi }{2}} \right)$ and $r \in \left( {0,1} \right)$ be the point of $\Gamma $ lying on the real axis (see Fig. \ref{geodes}).
Then the circle, $C$, passing through the points $z_1,z_2,r$ is given by
\[{x^2} + {y^2} + \frac{{1 - {r^2}}}{{r - \cos \theta }}x + \frac{{r\left( {r\cos \theta  - 1} \right)}}{{r - \cos \theta }} = 0\]
and has centre $K = \left( {\frac{{{r^2} - 1}}{{2\left( {r - \cos \theta } \right)}},0} \right)$, as illustrated in Fig. \ref{geodes}. Since the line passing through $K$ and $z_1$ is vertical to the tangent, $\varepsilon_1$, of the circle $C$ at $z_1$, we infer that
\[{\lambda _{{\varepsilon _1}}} = \frac{{2r\cos \theta  - 2{{\cos }^2}\theta  + 1 - {r^2}}}{{ - 2\left( {r - \cos \theta } \right)\sin \theta }},\]
where ${\lambda _{{\varepsilon _1}}}$ denotes the slope of ${\varepsilon _1}$. In addition, ${\varepsilon _1}$ is vertical to the tangent, $\varepsilon_2$, of $\partial \mathbb{D}$ at $z_1$ and thus
\[\frac{{\cos \theta }}{{\sin \theta }} \cdot \frac{{2r\cos \theta  - 2{{\cos }^2}\theta  + 1 - {r^2}}}{{2\left( {r - \cos \theta } \right)\sin \theta }} =  - 1\]
or 
\[r = \frac{{1 - \sin \theta }}{{\cos \theta }}.\]
Therefore,
\begin{equation}\label{metr}
{e^{ - {d_{\mathbb{D}}}\left( {0,\Gamma } \right)}} = \frac{{1 - r}}{{1 + r}} = \frac{{\cos \theta  + \sin \theta  - 1}}{{\cos \theta  - \sin \theta  + 1}}.
\end{equation}
Since the function 
\[f\left( \theta  \right) = \frac{{2\theta }}{\pi } \cdot  \frac{{\cos \theta  - \sin \theta  + 1}}{{\cos \theta  + \sin \theta  - 1}}\]
is decreasing on $\left( {0,\frac{\pi }{2}} \right)$ and
\[\mathop {\lim }\limits_{\theta  \to {0^ + }} f\left( \theta  \right) = \frac{4}{\pi },\;\mathop {\lim }\limits_{\theta  \to {{\frac{\pi }{2}}^ - }} f\left( \theta  \right) = 1,\]
we deduce that $1 \le f\left( \theta  \right) \le \frac{4}{\pi }$ for every $\theta  \in \left( {0,\frac{\pi }{2}} \right)$. This in conjunction with (\ref{metr}) and the fact that ${\omega _{\mathbb{D}}}\left( {0,\Gamma } \right) = \frac{{2\theta }}{\pi }$ (see \cite[p. 370]{Be}) gives the desired result.
\qed
\\\\
By the conformal invariance of harmonic measure, we can easily make the following computation.

\begin{lemma}\label{le2} Let $a,b \in \left( {0,1} \right)$. Then
\[{\omega _{\mathbb{D} \backslash \left[ {a,1} \right)}}\left( { - b,\partial {\mathbb{D}}} \right) = 1 - \frac{2}{\pi }\arctan \frac{1}{{\sqrt {{{\left( {\frac{{\left( {1 + a} \right)\left( {1 + b} \right)}}{{\left( {1 - a} \right)\left( {1 - b} \right)}}} \right)}^2} - 1} }}.\]
\end{lemma}

Hereinafter, let $\psi $ be a conformal map on $\mathbb{D}$ with $ \psi \left( 0 \right)=0$ and let ${F_\alpha } = \left\{ {z \in \mathbb{D}:\left| {\psi \left( z \right)} \right| = \alpha } \right\}$ and ${E_\alpha } = \left\{ {{e^{i\theta }}:\left| {\psi \left( {{e^{i\theta }}} \right)} \right| > \alpha } \right\}$ for $\alpha >0$. Moreover, set $d = \dist\left( {0,\partial \psi \left( \mathbb{D} \right)} \right)$ and let $N\left( \alpha \right) \in \mathbb{N} \cup \left\{ { + \infty } \right\}$ denote the number of components of $F_\alpha$ for $\alpha>0$.

\begin{lemma}\label{max} Let $F_\alpha ^i$ denote the components of $F_\alpha$, $i = 1,2, \ldots ,N\left( \alpha  \right)$. Then, for every $\alpha>0$, there exists a component $F_\alpha ^ *$ such that
\[{\omega _\mathbb{D}}\left( {0,F_\alpha ^ * } \right) = \max \left\{ {{\omega _\mathbb{D}}\left( {0,F_\alpha ^i} \right):i \in \left\{ {1,2, \ldots N\left( \alpha  \right)} \right\}} \right\}.\]	
\end{lemma}

\proof Fix an $\alpha>0$. Since the case $N\left( \alpha \right)<+\infty$ is trivial, suppose $N\left( \alpha \right)=+\infty$. Then the series 
\[\sum\limits_{i = 1}^{ + \infty } {{\omega _\mathbb{D}}\left( {0,F_\alpha ^i} \right)}  = {\omega _\mathbb{D}}\left( {0,{F_\alpha }} \right) \le 1\]
converges and hence
\[\mathop {\lim }\limits_{i \to  + \infty } {\omega _\mathbb{D}}\left( {0,F_\alpha ^i} \right) = 0.\]
This implies that $\exists {i_0} \in \mathbb{N}$ such that ${\omega _\mathbb{D}}\left( {0,F_\alpha ^i} \right) \le {\omega _\mathbb{D}}\left( {0,F_\alpha ^1} \right)$ for every $i \ge {i_0}$. So, setting ${\omega ^ * } = \max \left\{ {{\omega _\mathbb{D}}\left( {0,F_\alpha ^1} \right),{\omega _\mathbb{D}}\left( {0,F_\alpha ^2} \right), \ldots ,{\omega _\mathbb{D}}\left( {0,F_\alpha ^{{i_0} - 1}} \right)} \right\}$, we infer that there exists a component,  $F_\alpha ^ *$, of  $F_\alpha $ such that
\[{\omega _\mathbb{D}}\left( {0,F_\alpha ^ * } \right) ={\omega ^ * }= \max \left\{ {{\omega _\mathbb{D}}\left( {0,F_\alpha ^i} \right):i \in \left\{ {1,2, \ldots N\left( \alpha  \right)} \right\}} \right\}.\]
\qed

\begin{figure}[H] 
	\begin{center}
		\includegraphics[scale=0.52]{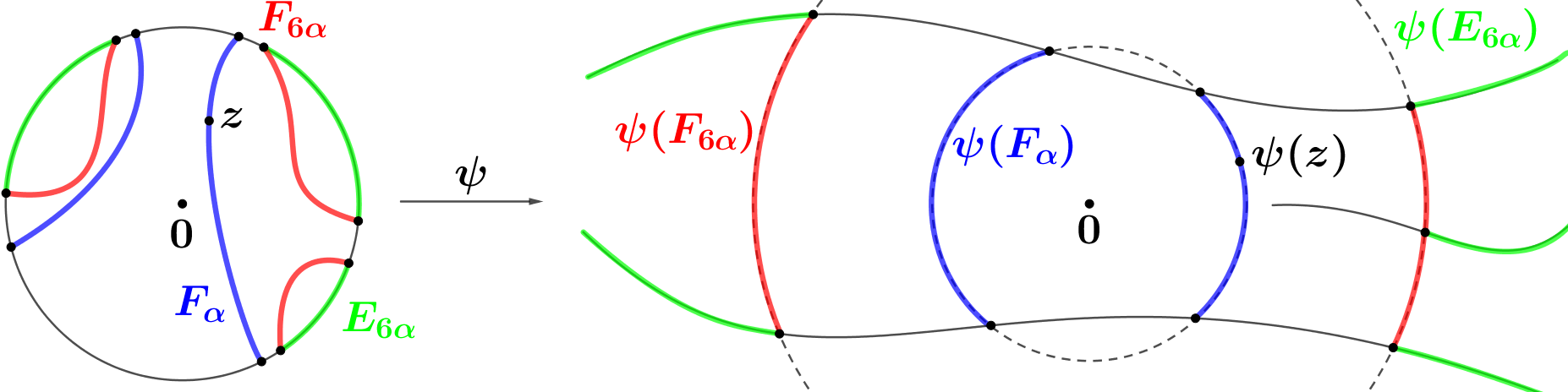}
		\caption{}
		\label{lem1}
	\end{center}
\end{figure}

\begin{lemma}\label{le3}
With the notation above, it is true that 
\[{\omega _\mathbb{D}}\left( {z,{E_{6\alpha }}} \right) \le \frac{1}{2},\;\forall z \in {F_\alpha },\;\forall \alpha  \ge 33d.\]
\end{lemma}

\begin{figure}[H] 
	\begin{center}
		\includegraphics[scale=0.55]{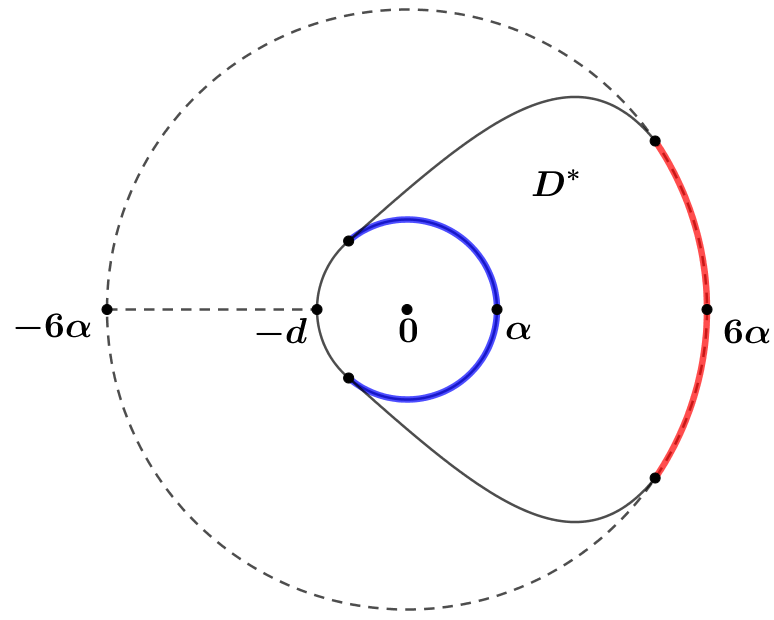}
		\caption{}
		\label{symm}
	\end{center}
\end{figure}

\proof Set $\psi \left( \mathbb{D} \right)=D$. If $z \in {F_\alpha }$ (see Fig. \ref{lem1}), then by Baernstein' s circular symmetrization (see \cite[Theorem 7]{Bae} and \cite[p. 665-669]{Hay}), Theorem \ref{pro} and the conformal invariance of harmonic mesaure, we infer that for every $\alpha  \ge 33d$,
\[{\omega _\mathbb{D}}\left( {z,{E_{6\alpha }}} \right) \le {\omega _\mathbb{D}}\left( {z,{F_{6\alpha }}} \right)={\omega _D}\left( {\psi \left( z \right),\psi \left( {{F_{6\alpha }}} \right)} \right) \le {\omega _{{D^ * }}}\left( {\alpha ,\partial {D^ * } \cap 6\alpha \partial \mathbb{D}} \right),\]
where $D^ *$ is the simply connected domain obtained by the circular symmetrization of $D \cap 6\alpha \mathbb{D}$ (see Fig. \ref{symm}). Applying Theorem \ref{pro}, the conformal invariance of harmonic mesaure and Lemma \ref{le2}, we have that for every $\alpha  \ge 33d$,
\begin{eqnarray}
{\omega _\mathbb{D}}\left( {z,{E_{6\alpha }}} \right) & \le& {\omega _{{D^ * }}}\left( {\alpha ,\partial {D^ * } \cap 6\alpha \partial \mathbb{D}} \right)  \le {\omega _{6\alpha \mathbb{D}\backslash \left( { - 6\alpha , - d} \right]}}\left( {\alpha ,6\alpha \partial \mathbb{D}} \right) = {\omega _{\mathbb{D}\backslash \left[ {\frac{d}{{6\alpha }},1} \right)}}\left( { - \frac{1}{6},\partial \mathbb{D}} \right) \nonumber \\ 
&=& 1 - \frac{2}{\pi }\arctan \frac{1}{{\sqrt {{{\left( {\frac{7}{5}\frac{{\left( {6a + d} \right)}}{{\left( {6a - d} \right)}}} \right)}^2} - 1} }} \le \frac{1}{2}, \nonumber 
\end{eqnarray}
where the last inequality comes from the fact that $\alpha  \ge \frac{{7 + 5\sqrt 2 }}{{30\sqrt 2  - 42}}d$. So, 
\[{\omega _\mathbb{D}}\left( {z,{E_{6\alpha }}} \right) \le \frac{1}{2},\;\forall z \in {F_\alpha },\;\forall \alpha  \ge 33d.\]
\qed

\begin{lemma}\label{le4} Let $c = \frac{{2 + \sqrt 2 }}{{2 - \sqrt 2 }}$ and $\alpha  > d$. Suppose that ${F_{c\alpha }^ * }$ is a component of ${F_{c\alpha }}$ such that
\[{\omega _\mathbb{D}}\left( {0,F_{c\alpha} ^ * } \right) = \max \left\{ {{\omega _\mathbb{D}}\left( {0,F_{c\alpha} ^i} \right):i \in \left\{ {1,2, \ldots N\left( \alpha  \right)} \right\}} \right\}\]  
and ${F_{\alpha }'}$ is the component of $F_\alpha $ such that ${F_{c\alpha }^ * }$ lies in the component of $\mathbb{D}\backslash F_{\alpha }'$ not containing the origin. If ${\Gamma_{\alpha }'}$ is the hyperbolic geodesic joining the endpoints of ${F_{\alpha }'}$ in $\mathbb{D}$, then
\[{\omega_\mathbb{D}}\left( {0,{F_{c\alpha }^ * }} \right) \le {\omega_\mathbb{D}}\left( {0,{\Gamma_{\alpha }'}} \right).\]
\end{lemma}

\begin{figure}[H] 
	\begin{center}
		\includegraphics[scale=0.5]{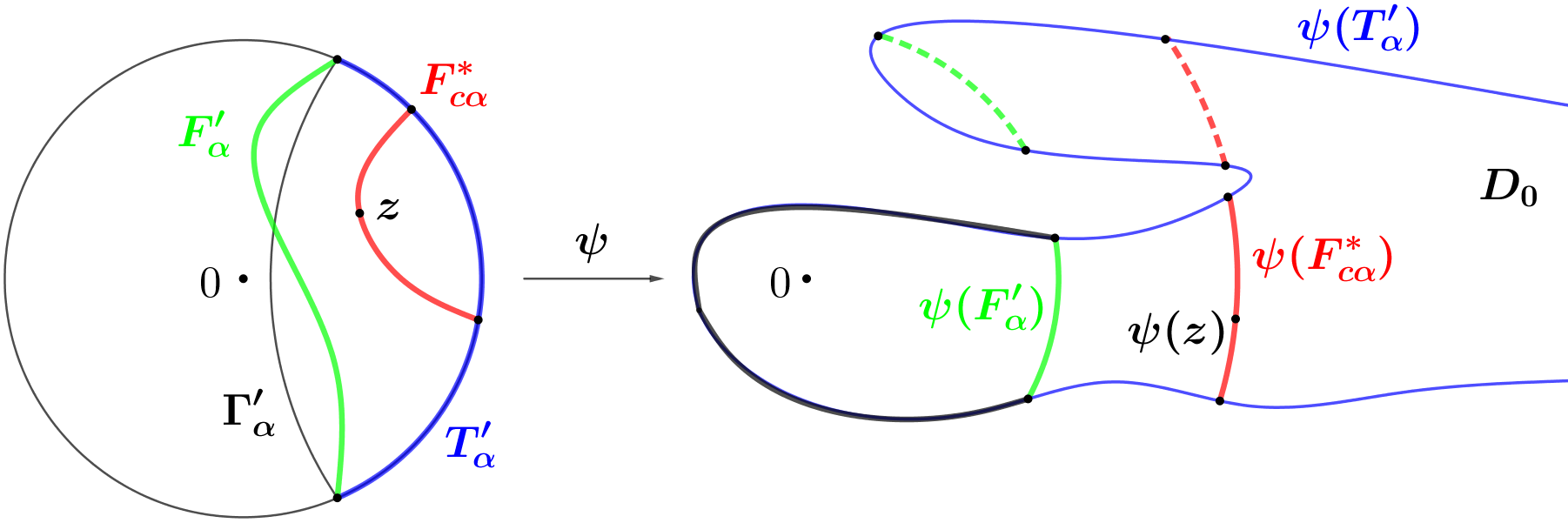}
		\caption{}
		\label{lem2}
	\end{center}
\end{figure}

\proof Lemma \ref{max} implies that a maximal component $ F_{c\alpha }^ *$ exists. Let $z \in F_{c\alpha }^ *$ and $\psi \left( \mathbb{D} \right) = D$. Let  ${T_{\alpha }'}$ be the arc of $\partial \mathbb{D}$ joining the endpoints of ${\Gamma_{\alpha }'}$ such that the interior of ${\Gamma _{\alpha }'} \cup {T_{\alpha }'}$ does not contain the origin (see Fig. \ref{lem2}). If $D_0$ is the component of $D\backslash \psi \left( {F_{\alpha }'} \right)$ containing $\psi \left( z \right)$, then 
\[{\omega _\mathbb{D}}\left( {z,{T_{\alpha }'} } \right) = {\omega _D}\left( {\psi \left( z \right),\psi \left( {T_{\alpha }'} \right)} \right) \ge {\omega _{{D_0}\backslash \alpha \overline {\mathbb{D}} }}\left( {\psi \left( z \right),\psi \left( {T_{\alpha }'} \right)}\backslash \alpha \overline {\mathbb{D}} \right).\]
Apply the conformal map $f\left( z \right) = \frac{\alpha }{z}$ which sends ${D_0}\backslash \alpha \overline {\mathbb{D}}$ onto an open set $W \subset \mathbb{D}$, the set $\left( {\mathbb{C} \backslash {\left( {D_0}\backslash \alpha \overline {\mathbb{D}} \right)}} \right)\backslash \left( {\alpha \overline {\mathbb{D}} } \right)$ into a set $A$ connecting $0$ to $\partial \mathbb{D}$ and $f\left( {\psi \left( z \right)} \right) = \frac{\alpha }{{\psi \left( z \right)}} \in \left( {\frac{1}{c}\partial \mathbb{D}} \right)$. Then, by the Beurling-Nevanlinna projection theorem, 
\[{\omega _{{D_0}\backslash \alpha \overline {\mathbb{D}} }}\left( {\psi \left( z \right),\psi \left( {T_{\alpha }'} \right)}\backslash \alpha \overline {\mathbb{D}} \right) = {\omega _W}\left( {\frac{\alpha }{{\psi \left( z \right)}},A} \right) \ge {\omega _\mathbb{D}}\left( {\frac{1}{c},\left( { - 1,0} \right]} \right) = \frac{2}{\pi }\arcsin \frac{{c - 1}}{{c + 1}} = \frac{1}{2}.\]
Therefore, 
\[{\omega _\mathbb{D}}\left( {z,{T_{\alpha }'} } \right) \ge \frac{1}{2},\;\forall z \in F_{c\alpha }^ *\]
which implies that $F_{c\alpha }^ *$ lies in the component of $\mathbb{D}\backslash \Gamma_{\alpha }'$ not containing the origin and hence
\[{\omega_\mathbb{D}}\left( {0,{F_{c\alpha }^ * }} \right) \le {\omega_\mathbb{D}}\left( {0,{\Gamma_{\alpha }'}} \right).\]
\qed

\section{Proofs}\label{section4}

\proof [Proof of Theorem \ref{hn}] The Beurling-Nevanlinna projection theorem implies that for every $\alpha>d$, 

\[{\omega _\mathbb{D}}\left( {0,{F_\alpha }} \right) \ge \frac{2}{\pi }{e^{ - {d_\mathbb{D}}\left( {0,{F_\alpha }} \right)}}\]
or equivalently
\[\frac{{\log {\omega _\mathbb{D}}{{\left( {0,{F_\alpha }} \right)}^{ - 1}}}}{{\log \alpha }} \le \frac{{\log \left( {{\pi  \mathord{\left/
 {\vphantom {\pi  2}} \right.
 \kern-\nulldelimiterspace} 2}} \right)}}{{\log \alpha }} + \frac{{{d_\mathbb{D}}\left( {0,{F_\alpha }} \right)}}{{\log \alpha }}.\]
(see \cite[p. 10]{Co}). By this and (\ref{sug}), we infer that
\[{\tt h}\left( \psi  \right) = \mathop {\liminf}\limits_{\alpha  \to  + \infty } \frac{{\log {\omega _\mathbb{D}}{{\left( {0,{F_\alpha }} \right)}^{ - 1}}}}{{\log \alpha }} \le \mathop {\liminf}\limits_{\alpha  \to  + \infty } \frac{{{d_\mathbb{D}}\left( {0,{F_\alpha }} \right)}}{{\log \alpha }}.\]
This in conjunction with the fact that 
\[\mathop {\liminf}\limits_{\alpha  \to  + \infty } \frac{{{d_\mathbb{\mathbb{D}}}\left( {0,{F_\alpha }} \right)}}{{\log \alpha }} \le {\tt h}\left( {\psi} \right),\]
which comes from Theorem \ref{Corra}, gives the desired result
\[{\tt h}\left( {\psi} \right) = \mathop {\liminf}\limits_{\alpha  \to  + \infty } \frac{{{d_\mathbb{D}}\left( {0,{F_\alpha }} \right)}}{{\log \alpha }}.\]
\qed

When the limits $\mathop {\lim }\limits_{\alpha  \to  + \infty } \frac{{\log {\omega _\mathbb{D}}{{\left( {0,{F_\alpha }} \right)}^{ - 1}}}}{{\log \alpha }}$ and $\mathop {\lim }\limits_{\alpha  \to  + \infty } \frac{{{d_\mathbb{D}}\left( {0,{F_\alpha }} \right)}}{{\log \alpha }}$ exist, we denote them by $L$ and $\mu $ respectively.

\proof[Proof of Corollary \ref{coth}] By Theorem \ref{hn} we obtain 
\begin{eqnarray}
\mu  =  + \infty & \Leftrightarrow& \mathop {\liminf}\limits_{\alpha  \to  + \infty } \frac{{{d_\mathbb{D}}\left( {0,{F_\alpha }} \right)}}{{\log \alpha }} =  + \infty  \Leftrightarrow \mathop {\liminf}\limits_{\alpha  \to  + \infty } \frac{{\log {\omega _\mathbb{D}}{{\left( {0,{F_\alpha }} \right)}^{ - 1}}}}{{\log \alpha }} =  + \infty \nonumber \\ 
&\Leftrightarrow& \mathop {\lim }\limits_{\alpha  \to  + \infty } \frac{{\log {\omega _\mathbb{D}}{{\left( {0,{F_\alpha }} \right)}^{ - 1}}}}{{\log \alpha }} =  + \infty \Leftrightarrow L =  + \infty. \nonumber 
\end{eqnarray}
\qed

\proof[Proof of Theorem \ref{thL}] If $\mu $ exists then Theorem \ref{hn} gives
\[\mathop {\liminf}\limits_{\alpha  \to  + \infty } \frac{{\log {\omega _\mathbb{D}}{{\left( {0,{F_\alpha }} \right)}^{ - 1}}}}{{\log \alpha }} = {\tt h}\left( {\psi} \right) = \mathop {\liminf}\limits_{\alpha  \to  + \infty } \frac{{{d_\mathbb{D}}\left( {0,{F_\alpha }} \right)}}{{\log \alpha }} = \mu. \]
By the Beurling-Nevanlinna projection theorem, for every $\alpha>d$,
\[\frac{{\log {\omega _\mathbb{D}}{{\left( {0,{F_\alpha }} \right)}^{ - 1}}}}{{\log \alpha }} \le \frac{{\log \left( {{\pi  \mathord{\left/
 {\vphantom {\pi  2}} \right.
 \kern-\nulldelimiterspace} 2}} \right)}}{{\log \alpha }} + \frac{{{d_\mathbb{D}}\left( {0,{F_\alpha }} \right)}}{{\log \alpha }}\]
and thus
\[\mathop {\lim \sup }\limits_{\alpha  \to  + \infty } \frac{{\log {\omega _\mathbb{D}}{{\left( {0,{F_\alpha }} \right)}^{ - 1}}}}{{\log \alpha }} \le \mathop {\lim \sup }\limits_{\alpha  \to  + \infty } \frac{{{d_\mathbb{D}}\left( {0,{F_\alpha }} \right)}}{{\log \alpha }} = \mu  = \mathop {\liminf}\limits_{\alpha  \to  + \infty } \frac{{\log {\omega _\mathbb{D}}{{\left( {0,{F_\alpha }} \right)}^{ - 1}}}}{{\log \alpha }}\]
which implies that
\[\mathop {\lim \sup }\limits_{\alpha  \to  + \infty } \frac{{\log {\omega _\mathbb{D}}{{\left( {0,{F_\alpha }} \right)}^{ - 1}}}}{{\log \alpha }} = \mathop {\liminf}\limits_{\alpha  \to  + \infty } \frac{{\log {\omega _\mathbb{D}}{{\left( {0,{F_\alpha }} \right)}^{ - 1}}}}{{\log \alpha }} = \mu. \]
So, $L$ exists and $L= \mu$.
\qed

\begin{figure}[H] 
	\begin{center}
		\includegraphics[scale=0.55]{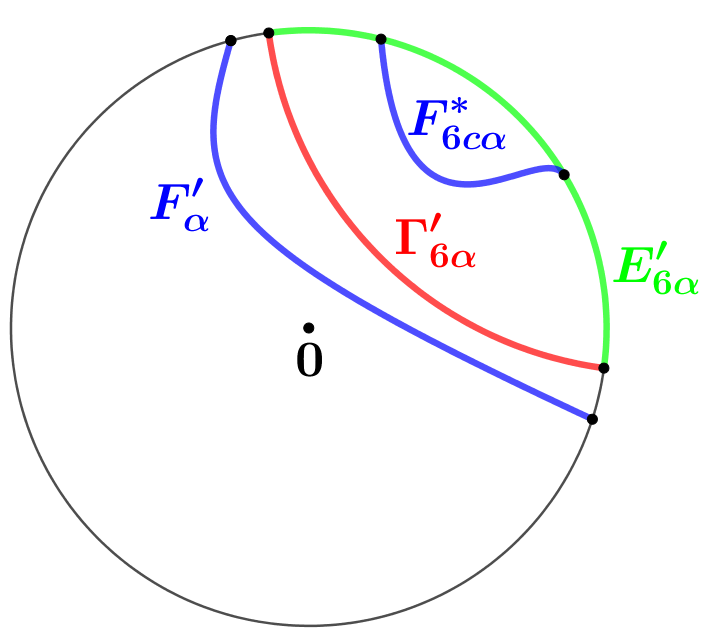}
		\caption{}
		\label{th}
	\end{center}
\end{figure}

\proof[Proof of Theorem \ref{cond}] If $L$ exists then Theorem \ref{hn} implies that
\begin{equation}\label{inf} 
\mathop {\liminf}\limits_{\alpha  \to  + \infty } \frac{{{d_\mathbb{D}}\left( {0,{F_\alpha }} \right)}}{{\log \alpha }} = \mathop {\liminf}\limits_{\alpha  \to  + \infty } \frac{{\log {\omega _\mathbb{D}}{{\left( {0,{F_\alpha }} \right)}^{ - 1}}}}{{\log \alpha }} = L.
\end{equation}
If $F_\alpha ^m$ denotes a component of $F_\alpha$ such that 
\[{d_\mathbb{D}}\left( {0,F_\alpha ^m} \right) = {d_\mathbb{D}}\left( {0,{F_\alpha }} \right),\]
then by the Beurling-Nevanlinna projection theorem, we get that for every $\alpha>d$,
\[{e^{ - {d_\mathbb{D}}\left( {0,{F_\alpha }} \right)}} = {e^{ - {d_\mathbb{D}}\left( {0,F_\alpha ^m} \right)}} \le \frac{\pi }{2}{\omega _\mathbb{D}}\left( {0,F_\alpha ^m} \right) \le \frac{\pi }{2}{\omega _\mathbb{D}}\left( {0,F_\alpha ^ * } \right)\]
or
\[\frac{{{d_\mathbb{D}}\left( {0,{F_\alpha }} \right)}}{{\log \alpha }} \ge \frac{{\log \left( {{2 \mathord{\left/
 {\vphantom {2 \pi }} \right.
 \kern-\nulldelimiterspace} \pi }} \right)}}{{\log \alpha }} + \frac{{\log {\omega _\mathbb{D}}{{\left( {0,F_\alpha ^ * } \right)}^{ - 1}}}}{{\log \alpha }}.\]
Thus
\begin{equation}\label{sup}
\mathop {\limsup }\limits_{\alpha  \to  + \infty } \frac{{{d_\mathbb{D}}\left( {0,{F_\alpha }} \right)}}{{\log \alpha }} \ge \mathop {\lim \sup }\limits_{\alpha  \to  + \infty } \frac{{\log {\omega _\mathbb{D}}{{\left( {0,F_\alpha ^ * } \right)}^{ - 1}}}}{{\log \alpha }}.
\end{equation}
If $c = \frac{{2 + \sqrt 2 }}{{2 - \sqrt 2 }}$, let ${F_{\alpha }'}$ be the component of $F_\alpha $ such that ${F_{6c\alpha }^ * }$ lies in the component of $\mathbb{D}\backslash F_{\alpha }'$ not containing the origin (see Fig. \ref{th}). Also, let $E _{6\alpha }'$ be the arc of $E _{6\alpha }$ such that $E _{6\alpha }' \cap \overline F _{6c\alpha }^ *  \ne \emptyset$ and $\Gamma _{6\alpha }'$ be the hyperbolic geodesic joining the endpoints of $E _{6\alpha }'$. By Lemma \ref{le3} we have that 
\[{\omega _\mathbb{D}}\left( {z,{E_{6\alpha }}} \right) \le \frac{1}{2},\;\forall z \in {F_\alpha },\;\forall \alpha  \ge 33d\]
and thus
\[{\omega _\mathbb{D}}\left( {z,{E_{6\alpha }'}} \right) \le \frac{1}{2},\;\forall z \in {F_{\alpha}' },\;\forall \alpha  \ge 33d.\]
This implies that $F_{\alpha}'$ lies in the component of $\mathbb{D}\backslash \Gamma_{6\alpha }'$ containing the origin and hence
\[{d_\mathbb{D}}\left( {0,{F_{\alpha}' }} \right) \le {d_\mathbb{D}}\left( {0,{\Gamma_{6\alpha}' }} \right).\]
This and Lemma \ref{geodesic} give that for every $\alpha  \ge 33d$,
\begin{equation}\label{1} {e^{ - {d_\mathbb{D}}\left( {0,{F_\alpha }} \right)}} \ge {e^{ - {d_\mathbb{D}}\left( {0,{F_{\alpha}' }} \right)}} \ge {e^{ - {d_\mathbb{D}}\left( {0,{\Gamma_{6\alpha}' }} \right)}} \ge \frac{\pi }{4}{\omega _\mathbb{D}}\left( {0,{\Gamma_{6\alpha}'}} \right).
\end{equation}
By Lemma \ref{le4}, we get 
\begin{equation}\label{2} {\omega _\mathbb{D}}\left( {0,{\Gamma _{6\alpha }'}} \right) \ge {\omega _\mathbb{D}}\left( {0,F_{6c\alpha }^ * } \right).
\end{equation}
Combining the relations (\ref{1}) and (\ref{2}), we infer that for every $\alpha  \ge 33d$,
\[{e^{ - {d_\mathbb{D}}\left( {0,{F_\alpha }} \right)}} \ge \frac{\pi }{4}{\omega _\mathbb{D}}\left( {0,F_{6c\alpha }^ * } \right),\] 
or equivalently
\[\frac{{{d_\mathbb{D}}\left( {0,{F_\alpha }} \right)}}{{\log \alpha }} \le \frac{{\log \left( {{4 \mathord{\left/
 {\vphantom {4 \pi }} \right.
 \kern-\nulldelimiterspace} \pi }} \right)}}{{\log \alpha }} + \frac{{\log {\omega _\mathbb{D}}{{\left( {0,F_{6c\alpha }^ * } \right)}^{ - 1}}}}{{\log \alpha }}.\]
Therefore,
\begin{eqnarray}
\mathop {\lim \sup }\limits_{\alpha  \to  + \infty } \frac{{{d_\mathbb{D}}\left( {0,{F_\alpha }} \right)}}{{\log \alpha }} &\le& \mathop {\lim \sup }\limits_{\alpha  \to  + \infty } \frac{{\log {\omega _\mathbb{D}}{{\left( {0,F_{6c\alpha }^ * } \right)}^{ - 1}}}}{{\log \alpha }} = \mathop {\lim \sup }\limits_{\alpha  \to  + \infty } \left( {\frac{{\log {\omega _\mathbb{D}}{{\left( {0,F_{6c\alpha }^ * } \right)}^{ - 1}}}}{{\log \left( {6c\alpha } \right)}}\frac{{\log \left( {6c\alpha } \right)}}{{\log \alpha }}} \right) \nonumber \\
&=&\mathop {\lim \sup }\limits_{\alpha  \to  + \infty } \frac{{\log {\omega _\mathbb{D}}{{\left( {0,F_{6c\alpha }^ * } \right)}^{ - 1}}}}{{\log \left( {6c\alpha } \right)}} = \mathop {\lim \sup }\limits_{\alpha  \to  + \infty } \frac{{\log {\omega _\mathbb{D}}{{\left( {0,F_\alpha ^ * } \right)}^{ - 1}}}}{{\log \alpha }}. \nonumber
\end{eqnarray}
This in conjunction with (\ref{sup}) gives
\begin{equation}\label{ssup} \mathop {\limsup }\limits_{\alpha  \to  + \infty } \frac{{{d_\mathbb{D}}\left( {0,{F_\alpha }} \right)}}{{\log \alpha }} = \mathop {\lim \sup }\limits_{\alpha  \to  + \infty } \frac{{\log {\omega _\mathbb{D}}{{\left( {0,F_\alpha ^ * } \right)}^{ - 1}}}}{{\log \alpha }}.
\end{equation}
By relations (\ref{inf}) and (\ref{ssup}), we conclude that $\mu$ exists if and only if 
\[\mathop {\limsup }\limits_{\alpha  \to  + \infty } \frac{{{d_\mathbb{D}}\left( {0,{F_\alpha }} \right)}}{{\log \alpha }} = L \Leftrightarrow  \mathop {\lim \sup }\limits_{\alpha  \to  + \infty } \frac{{\log {\omega _\mathbb{D}}{{\left( {0,F_\alpha ^ * } \right)}^{ - 1}}}}{{\log \alpha }}=L\]
and if $\mu$ exists then $\mu=L$.
\qed

\proof[Proof of Corollary \ref{cor}] Obviously, for every $\alpha>0$,
\[\frac{1}{{N\left( \alpha  \right)}}{\omega _\mathbb{D}}\left( {0,{F_\alpha }} \right) \le {\omega _\mathbb{D}}\left( {0,F_\alpha ^ * } \right) \le {\omega _\mathbb{D}}\left( {0,{F_\alpha }} \right)\]
or
\[{\omega _\mathbb{D}}{\left( {0,{F_\alpha }} \right)^{ - 1}} \le {\omega _\mathbb{D}}{\left( {0,F_\alpha ^ * } \right)^{ - 1}} \le N\left( \alpha  \right){\omega _\mathbb{D}}{\left( {0,{F_\alpha }} \right)^{ - 1}}\]
or
\[\frac{{\log {\omega _\mathbb{D}}{{\left( {0,{F_\alpha }} \right)}^{ - 1}}}}{{\log \alpha }} \le \frac{{\log {\omega _\mathbb{D}}{{\left( {0,F_\alpha ^ * } \right)}^{ - 1}}}}{{\log \alpha }} \le \frac{{\log N\left( \alpha  \right)}}{{\log \alpha }} + \frac{{\log {\omega _\mathbb{D}}{{\left( {0,{F_\alpha }} \right)}^{ - 1}}}}{{\log \alpha }}.\]
Since $L$ exists and $\mathop {\lim }\limits_{\alpha  \to  + \infty } \frac{\log{N\left( \alpha  \right)}}{{\log \alpha }} = 0$, the above inequalities give that
\[\mathop {\lim }\limits_{\alpha  \to  + \infty } \frac{{\log {\omega _\mathbb{D}}{{\left( {0,F_\alpha ^ * } \right)}^{ - 1}}}}{{\log \alpha }} = L\]
and thus Theorem \ref{cond} implies that $\mu $ exists and $\mu=L $.
\qed

\section{Example}\label{section5}

\begin{figure}[H]
	\begin{minipage}{0.4\textwidth}
		\begin{center}
			\includegraphics[scale=0.5]{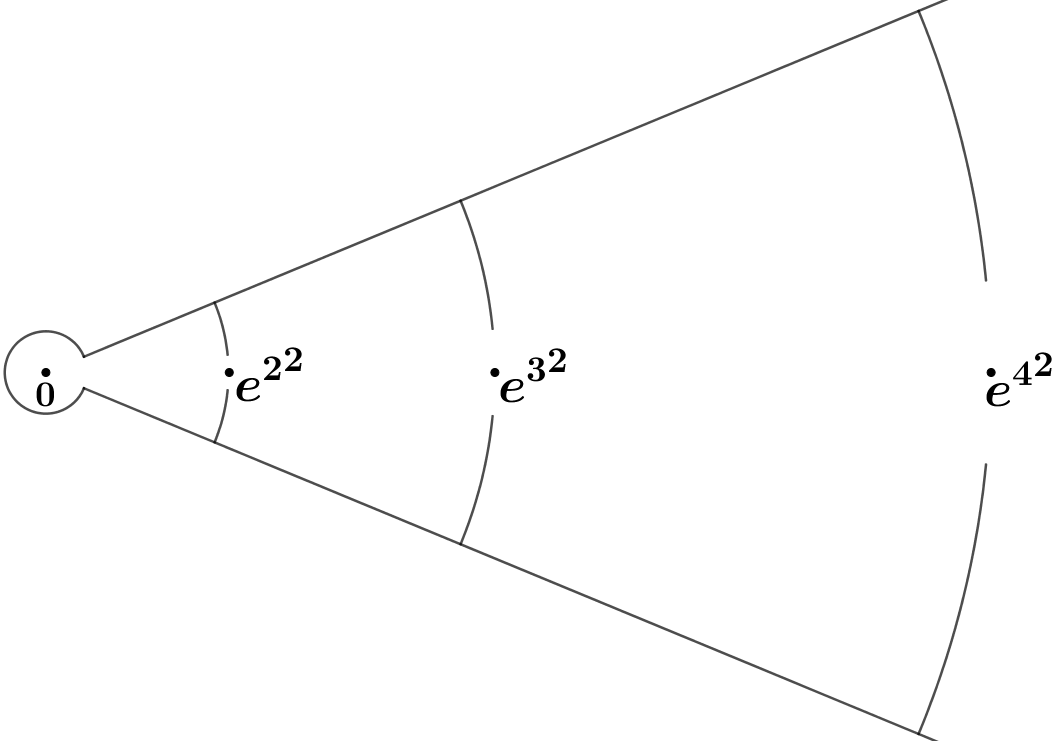}
			\caption{The simply connected domain $D$.}
			\label{t2}
		\end{center}
	\end{minipage}\hfill
	\begin{minipage}{0.5\textwidth}
		\begin{center}
			\includegraphics[scale=0.5]{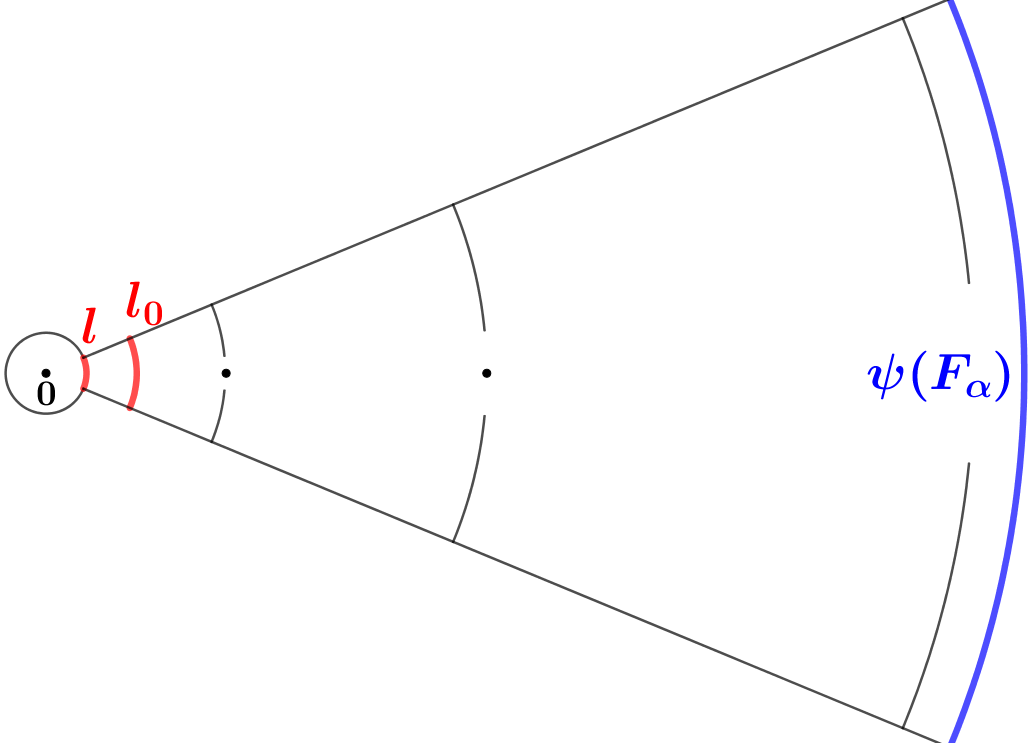}
			\vspace*{0.2cm}
			\caption{The crosscuts $l,l_0$ and $\psi \left( {{F_\alpha }} \right)$ in $D$.}
			\label{t3}
		\end{center}
	\end{minipage}
\end{figure}

\begin{figure}[H]
	\begin{center}
		\includegraphics[scale=0.55]{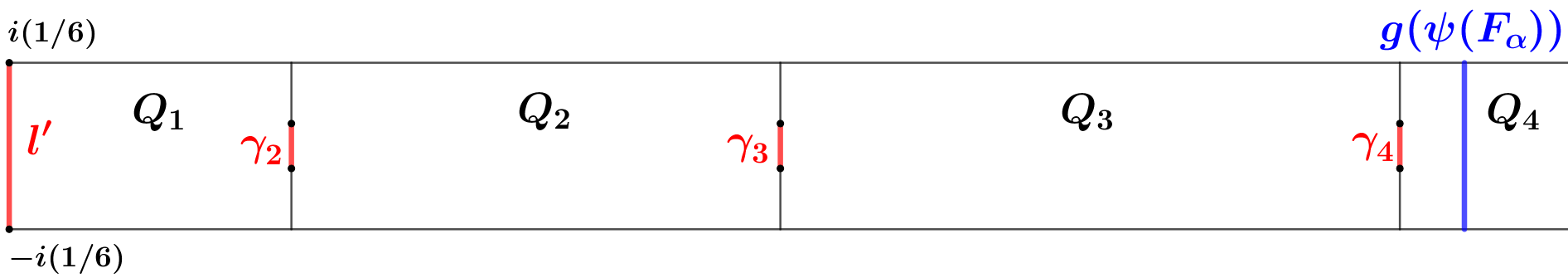}
		\caption{The crosscuts $\gamma _j$ and the quadrilaterals $Q_j$ in ${D_0^\prime}$.}
		\label{rect}
	\end{center}
\end{figure}

Trying to find a conformal map $\psi$ on $\mathbb{D}$ such that $\psi \left( 0 \right) = 0$, ${\tt h} \left( \psi  \right) = \mu  = L <  + \infty $ and $\psi  \in {H^\mu }\left( \mathbb{D} \right)$, we had to deal with the following issues:
\begin{enumerate}
\item $\mathop {\lim }\limits_{\alpha  \to  + \infty } \frac{\log{N\left( \alpha  \right)}}{{\log \alpha }} = 0$ so as to ensure the existence of both $\mu$ and $L$ (see Theorem \ref{thL} and Corollary \ref{cor}),
\item Find exactly the number ${\tt h} \left( \psi  \right)$,
\item $ \int_0^{ + \infty } {{\alpha ^{\mu - 1}}{\omega _{\mathbb{D}}}\left( {0,{F_\alpha }} \right)d\alpha }  <  + \infty$ so that $\psi  \in {H^{\mu}}\left( \mathbb{D} \right)$ (see the relation (\ref{isod})).
\end{enumerate}
So, considering the simply connected domain $D$ of Fig. \ref{t2} and the corresponding Riemann map $\psi$ from $\mathbb{D}$ onto $D$ with $\psi \left( 0 \right) = 0$, we obtain:
\begin{enumerate}[\rm (i)]
	\item $N\left( \alpha  \right) = 1$ for every $\alpha>0$; so (1) is satisfied.
	\item The evaluation of ${\tt h} \left( \psi  \right)$ by estimating ${\omega _{\mathbb{D}}}\left( {0,{F_\alpha }} \right)$ with the aid of extremal length (see Theorem \ref{mar}) which can be estimated in a domain of the form illustrated in Fig. \ref{rect}, by  applying the serial rule and the domain decomposition method; so (2) is satisfied. 
	\item $ \int_0^{ + \infty } {{\alpha ^{\mu - 1}}{\omega _{\mathbb{D}}}\left( {0,{F_\alpha }} \right)d\alpha }  <  + \infty$ because of the circular arcs of $\partial D$ and because of the choice of the sequence $\left\{ {{e^{{n^2}}}} \right\}$ (see Fig. \ref{t2}) which we made after some trials; and thus (3) is satisfied.
\end{enumerate} 

\begin{example} There exists a conformal map $\psi$ on $\mathbb{D}$ such that $\mu$ exists and $\psi  \in {H^\mu {\left( \mathbb{D} \right)} }$.
\end{example}

\begin{figure}[H]
	\begin{center}
	\includegraphics[scale=0.55]{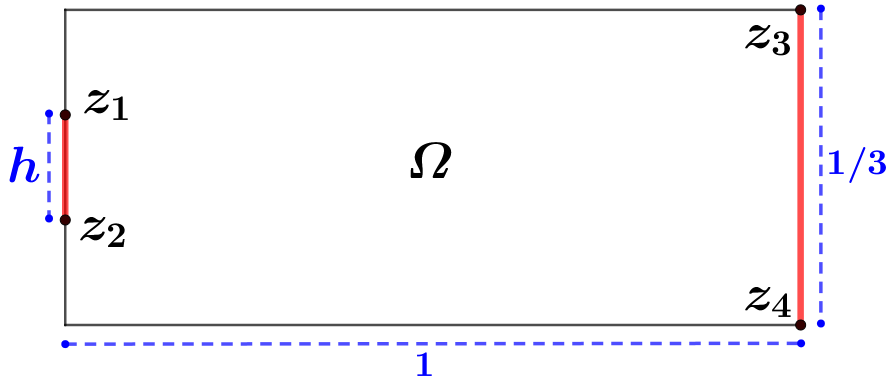}
	\caption{The quadrilateral  ${Q^ * } = \left\{ {\Omega ;{z_1},{z_2},{z_3},{z_4}} \right\}$.}
	\label{r1}
\end{center}
\end{figure}

\proof \textbf{Step 1:} Let $D$ be the simply connected domain of Fig. \ref{t2}, namely
\[D =\mathbb{D} \cup \left\{ {z \in \mathbb{C} :\left| {\Arg{z}} \right| < \frac{1}{6}} \right\}\backslash \bigcup\limits_{n = 2}^{ + \infty } {\left\{ {z \in  {e^{{n^2}}}\partial {\mathbb{D}} :\frac{h}{2} \le \left| {\Arg{z}} \right| < \frac{1}{6}} \right\}},\]
where $h$ is a positive constant small enough so that if $m\left( {{Q^ * }} \right)$ is the module of the quadrilateral ${Q^ * } = \left\{ {\Omega ;{z_1},{z_2},{z_3},{z_4}} \right\}$ illustrated in Fig. \ref{r1}, then $m\left( {{Q^ * }} \right) > 9$. The Riemann Mapping Theorem implies that there exists a conformal map  $\psi$ from $\mathbb{D}$ onto $D$ such that $\psi \left( 0 \right) = 0$.

\begin{figure}[H]
		\begin{center}
			\includegraphics[scale=0.54]{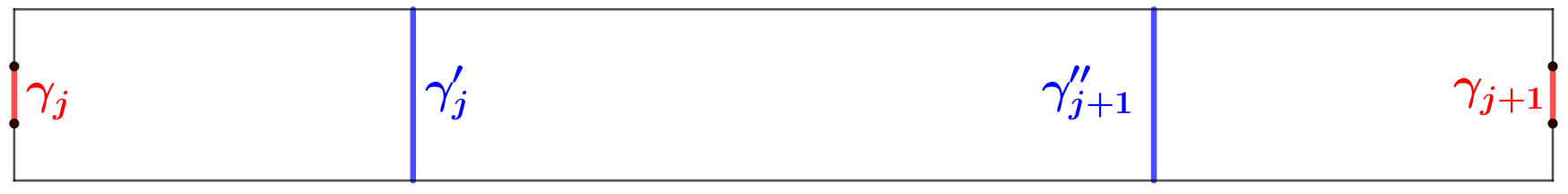}
			\caption{The crosscuts ${\gamma _j^\prime},{\gamma _{j + 1}^{\prime \prime }}$ in $Q_j$.}
			\label{r2}
		\end{center}
\end{figure}

\begin{figure}[H]
	\begin{center}
		\includegraphics[scale=0.5]{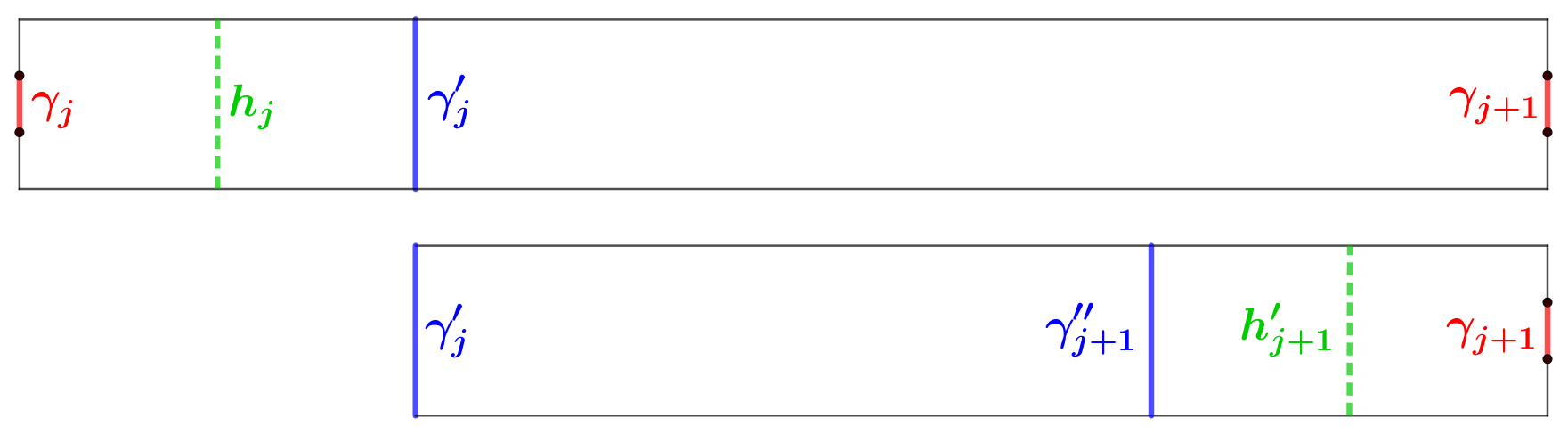}
		\caption{The auxilary crosscuts ${h_j},{h_{j + 1}^\prime}$.}
		\label{r3}
	\end{center}
\end{figure}

\textbf{Step 2:} Fix a real number $\alpha > {e^{{3^2}}}$. Then there exists a fixed number $n \in \mathbb{N}$ such that 
\begin{equation}\label{na} {e^{{n^2}}} \le \alpha  < {e^{{{\left( {n + 1} \right)}^2}}} \Leftrightarrow n \le \sqrt {\log \alpha }  < n + 1.
\end{equation}
Applying Theorems \ref{el} and \ref{mar}, we have 
\begin{equation}\label{sx1} {\omega _\mathbb{D}}\left( {0,{F_\alpha }} \right) = {\omega _D}\left( {0,\psi \left( {{F_\alpha }} \right)} \right) \le \frac{8}{\pi }{e^{ - \pi {\lambda _D}\left( {\left( { - 1,0} \right],\psi \left( {{F_\alpha }} \right)} \right)}} \le \frac{8}{\pi }{e^{ - \pi {\lambda _{{D_0}}}\left( {l,\psi \left( {{F_\alpha }} \right)} \right)}},
\end{equation}
where ${D_0} = D\backslash \overline {\mathbb{D}}$ and $l = \partial {\mathbb{D}} \cap D$ (see Fig. \ref{t3}). Set for $j = 2,3, \ldots,n+1$,
\[{\gamma _j} = \left\{ {{j^2} + iy:\left| y \right| \le \frac{h}{2}} \right\}.\]
Applying the conformal map $g\left( z \right) = \LOG \left( z \right)$ on $D_0$ and setting $g\left( {{D_0}} \right) = {D_0^\prime} $ and $g\left( l \right) = l'=\left\{ {iy:\left| y \right| \le \frac{1}{6}} \right\}$, we get by the conformal invariance of extremal length and Theorem \ref{el} that
\[{\lambda _{{D_0}}}\left( {l,\psi \left( {{F_\alpha }} \right)} \right) = {\lambda _{{D_0^\prime} }}\left( {l',g\left( {\psi \left( {{F_\alpha }} \right)} \right)} \right) \ge {\lambda _{{D_0^\prime} }}\left( {l',{\gamma _n}} \right).\]
This and (\ref{sx1}) give
\begin{equation}\label{sx2} {\omega _\mathbb{D}}\left( {0,{F_\alpha }} \right) \le \frac{8}{\pi }{e^{ - \pi {\lambda _{{D_0^\prime} }}\left( {l',{\gamma _n}} \right)}}.
\end{equation}
Taking the crosscuts ${\gamma _2},{\gamma _3}, \ldots, {\gamma _{n+1}}$ of ${D_0^\prime}$ and setting
\[m\left( {{Q_1}} \right) = {\lambda _{{D_0^\prime} }}\left( {l',{\gamma _2}} \right),\,m\left( {{Q_2}} \right) = {\lambda _{{D_0^\prime} }}\left( {{\gamma _2},{\gamma _3}} \right), \ldots,\, m\left( {{Q_{n}}} \right) = {\lambda _{{D_0^\prime} }}\left( {{\gamma _{n }},{\gamma _{n+1}}} \right)\]
as illustrated in Fig. \ref{rect}, the serial rule implies that
\begin{equation}\label{sx3} {\lambda _{{D_0^\prime} }}\left( {l',{\gamma _n}} \right) \ge m\left( {{Q_1}} \right) + m\left( {{Q_2}} \right) +  \ldots  + m\left( {{Q_{n - 1}}} \right) \ge m\left( {{Q_2}} \right) +  \ldots  + m\left( {{Q_{n - 1}}} \right).
\end{equation}
In every $Q_j$, for $j=2,3,\ldots,n$, we take the crosscuts
\[{\gamma _j^\prime}  = \left\{ {\left( {{j^2} + 1} \right) + iy:\left| y \right| \le \frac{1}{6}} \right\},\;{\gamma _{j + 1}^{\prime \prime }} = \left\{ {{{\left( {j + 1} \right)}^2} - 1 + iy:\left| y \right| \le \frac{1}{6}} \right\}\]
(see Fig.\ \ref{r2}) so that applying the serial rule,
\[m\left( {{Q_j}} \right) \ge 2m\left( {{Q^ * }} \right) + {\lambda _{{Q_j}}}\left( {{\gamma _j^\prime} ,{\gamma _{j + 1}^{\prime \prime }}} \right) = 2m\left( {{Q^ * }} \right) + 3\left( {2j - 1} \right).\]
Adding for $j=2,3,\ldots,n-1$, we get
\begin{eqnarray}
m\left( {{Q_2}} \right) +  \ldots  + m\left( {{Q_{n - 1}}} \right) &\ge& 2\left( {n - 2} \right)m\left( {{Q^ * }} \right) + 3\left( {3 + 5 + 7 +  \ldots \left( {2n - 3} \right)} \right) \nonumber \\
&=& 2\left( {n - 2} \right)m\left( {{Q^ * }} \right) + 3n{\left( {n - 2} \right)} \nonumber \\
&\ge& 2\left( {n - 2} \right)m\left( {{Q^ * }} \right) + 3{\left( {n - 2} \right)}^2 \nonumber \\
&\ge& 2\left( {\sqrt {\log \alpha }  - 3} \right)m\left( {{Q^ * }} \right) + 3{\left( {\sqrt {\log \alpha }  - 3} \right)^2}, \label{sx5}
\end{eqnarray}
where the last inequality comes from (\ref{na}).  Combining the relations (\ref{sx2}), (\ref{sx3}) and (\ref{sx5}), we infer that
\begin{equation}\label{st1} {\omega _\mathbb{D}}\left( {0,{F_\alpha }} \right) \le \frac{8}{\pi }{e^{ - 2\pi \left( {\sqrt {\log \alpha }  - 3} \right)m\left( {{Q^ * }} \right) - 3\pi {{\left( {\sqrt {\log \alpha }  - 3} \right)}^2}}}
\end{equation}
or equivalently
\begin{eqnarray} \frac{{\log {\omega _\mathbb{D}}{{\left( {0,{F_\alpha }} \right)}^{ - 1}}}}{{\log \alpha }} &\ge& \frac{{\log \left( {{\pi  \mathord{\left/
					{\vphantom {\pi  8}} \right.
					\kern-\nulldelimiterspace} 8}} \right)}}{{\log \alpha }} + \frac{{2\pi \left( {\sqrt {\log \alpha }  - 3} \right)m\left( {{Q^ * }} \right) + 3\pi {{\left( {\sqrt {\log \alpha }  - 3} \right)}^2}}}{{\log \alpha }}  \nonumber \\
&=& 3\pi  + \frac{{2\pi \left( {m\left( {{Q^ * }} \right) - 9} \right)\sqrt {\log \alpha }  + \log \left( {{\pi  \mathord{\left/
					{\vphantom {\pi  8}} \right.
					\kern-\nulldelimiterspace} 8}} \right) + 27\pi  - 6\pi m\left( {{Q^ * }} \right)}}{{\log \alpha }}. \label{kla}			
\end{eqnarray}
So, taking limits as $\alpha  \to  + \infty$,
\begin{equation}\label{liminf}\mathop {\liminf}\limits_{\alpha  \to  + \infty } \frac{{\log {\omega _\mathbb{D}}{{\left( {0,{F_\alpha }} \right)}^{ - 1}}}}{{\log \alpha }} \ge 3\pi.
\end{equation}
\textbf{Step 3:} On the other hand, by Theorem \ref{mar}, we have
\begin{equation}\label{sx11} {\omega _\mathbb{D}}\left( {0,{F_\alpha }} \right) = {\omega _D}\left( {0,\psi \left( {{F_\alpha }} \right)} \right) \ge {e^{ - \pi {\lambda _D}\left( {\left( { - 1,0} \right],\psi \left( {{F_\alpha }} \right)} \right)}}.
\end{equation}
Take the crosscut $l_0 = e\partial \mathbb{D} \cap D$ (see Fig. \ref{t3}). Then 
\[{\lambda _{{D_0}}}\left( {l,{l_0}} \right) = {\lambda _{{D_0^\prime} }}\left( {g\left( l \right),g\left( {{l_0}} \right)} \right) = 3\]
and thus Theorem \ref{ddm1} implies that
\begin{equation}\label{sx12}{\lambda _D}\left( {\left( { - 1,0} \right],\psi \left( {{F_\alpha }} \right)} \right) \le {C_0} + {\lambda _{{D_0}}}\left( {l,\psi \left( {{F_\alpha }} \right)} \right) + 2.71{e^{ - 3\pi }},
\end{equation}
where ${C_0}: = {\lambda _D}\left( {\left( { - 1,0} \right],{l_0}} \right)-3$. By Theorem \ref{el}, we take
\[{\lambda _{{D_0}}}\left( {l,\psi \left( {{F_\alpha }} \right)} \right) = {\lambda _{{D_0^\prime}}}\left( {l',g\left( {\psi \left( {{F_\alpha }} \right)} \right)} \right) \le {\lambda _{{D_0^\prime} }}\left( {l',{\gamma _{n+1}}} \right)\]
which gives with (\ref{sx11}) and (\ref{sx12}) that
\begin{equation}\label{sx22}
{\omega _\mathbb{D}}\left( {0,{F_\alpha }} \right) \ge{e^{ - K}}{e^{ - \pi {\lambda _{{D_0^\prime}}}\left( {l',{\gamma _{n + 1}}} \right)}},
\end{equation}
where $K: = {C_0}\pi  + 2.71{e^{ - 3\pi }}\pi $.
Considering the crosscuts ${\gamma _2},{\gamma _3}, \ldots, {\gamma _{n + 1}}$ of ${D_0^\prime}$ and applying successively Theorem \ref{ddm} by using every time the auxilary crosscuts ${\gamma _j^\prime} $ and ${\gamma _j^{\prime \prime }}$, we obtain 
\begin{eqnarray}
{\lambda _{{D_0^\prime} }}\left( {l',{\gamma _{n + 1}}} \right) &\le& m\left( {{Q_1}} \right) + m\left( {{{\left( {{Q_1}} \right)}^c}} \right) + 5.26{e^{ - 2\pi m\left( {{Q^ * }} \right)}}      \nonumber \\
m\left( {{{\left( {{Q_1}} \right)}^c}} \right) &\le& m\left( {{Q_2}} \right) + m\left( {{{\left( {{Q_2}} \right)}^c}} \right) + 5.26{e^{ - 2\pi m\left( {{Q^ * }} \right)}} \nonumber \\
\vdots \nonumber \\
m\left( {{{\left( {{Q_{n - 2}}} \right)}^c}} \right) &\le& m\left( {{Q_{n - 1}}} \right) + m\left( {{Q_n}} \right) + 5.26{e^{ - 2\pi m\left( {{Q^ * }} \right)}},  \nonumber
\end{eqnarray}
where $m\left( {{{\left( {{Q_j}} \right)}^c}} \right):={\lambda _{{D_0^\prime} }}\left( {{\gamma _{j + 1}},{\gamma _{n + 1}}} \right)$ for $j = 1,2, \ldots ,n - 2$. Adding the inequalities above, we deduce that
\begin{equation}\label{sx23}{\lambda _{{D_0^\prime} }}\left( {l',{\gamma _{n + 1}}} \right) \le m\left( {{Q_1}} \right) + m\left( {{Q_2}} \right) +  \ldots  + m\left( {{Q_n}} \right) + 5.26{e^{ - 2\pi m\left( {{Q^ * }} \right)}}\left( {n - 1} \right).
\end{equation}
Now set for $j=2,3,\ldots,n$,
\[{h_j} = \left\{ {\left( {{j^2} + \frac{1}{2}} \right) + iy:\left| y \right| \le \frac{1}{6}} \right\},\;{h_{j + 1}^\prime } = \left\{ {{{\left( {j + 1} \right)}^2} - \frac{1}{2} + iy:\left| y \right| \le \frac{1}{6}} \right\}.\]
In every $Q_j$, for $j=2,3,\ldots,n$, we take the crosscut ${\gamma _j^\prime} $ and the auxilary crosscut ${h_j}$ (see Fig. \ref{r3}). Since ${\lambda _{{D_0^\prime} }}\left( {{h_j},{\gamma _j^\prime} } \right) = \frac{3}{2}$, by applying Theorem \ref{ddm} we take
\[m\left( {{Q_j}} \right) \le m\left( {{Q^*}} \right) + {\lambda _{{D_0}^\prime }}\left( {{\gamma _j^\prime },{\gamma _{j + 1}}} \right) + 5.26{e^{ - 3\pi }}.\]
Then considering the crosscut ${\gamma _{j + 1}^{\prime \prime }}$ and the auxilary crosscut ${h_{j + 1}^\prime} $  (see Fig. \ref{r3}), we have again by Theorem \ref{ddm} that
\[{\lambda _{{D_0^\prime} }}\left( {{\gamma _j^\prime} ,{\gamma _{j + 1}}} \right) \le m\left( {{Q^*}} \right) + 3\left( {2j - 1} \right) + 5.26{e^{ - 3\pi }},\]
where ${\lambda _{{D_0^\prime} }}\left( {{\gamma _j^\prime} ,{\gamma _{j + 1}^{\prime \prime }}} \right) = 3\left( {2j - 1} \right)$. Combining the inequalities above, we finally get
\[m\left( {{Q_j}} \right) \le 2m\left( {{Q^ * }} \right) + 10.52{e^{ - 3\pi }} + 3\left( {2j - 1} \right).\]
This in conjunction with (\ref{sx23}) gives 
\begin{eqnarray}
{\lambda _{{D_0^\prime} }}\left( {l',{\gamma _{n + 1}}} \right) &\le& m\left( {{Q_1}} \right) + \left( {2m\left( {{Q^ * }} \right) + 10.52{e^{ - 3\pi }} + 5.26{e^{ - 2\pi m\left( {Q^*} \right)}}} \right)\left( {n - 1} \right) + 3\sum\limits_{j = 2}^n {\left( {2j - 1} \right)}       \nonumber \\
&=& m\left( {{Q_1}} \right) + \left( {2m\left( {{Q^ * }} \right) + 10.52{e^{ - 3\pi }} + 5.26{e^{ - 2\pi m\left( {Q^*} \right)}}} \right)\left( {n - 1} \right) + 3\left( {n - 1} \right)\left( {n + 1} \right)  \nonumber \\
&=& m\left( {{Q_1}} \right)-3 + \left( {2m\left( {{Q^ * }} \right) + 10.52{e^{ - 3\pi }} + 5.26{e^{ - 2\pi m\left( {Q^*} \right)}}} \right)\left( {n - 1} \right) + 3n^2  \nonumber \\
&\le& m\left( {{Q_1}} \right)-3 + \left( {2m\left( {{Q^ * }} \right) + 10.52{e^{ - 3\pi }} + 5.26{e^{ - 2\pi m\left( {Q^*} \right)}}} \right)\left( {\sqrt {\log \alpha }  - 1} \right) + 3{\log \alpha }, \nonumber
\end{eqnarray}
where the last inequality comes from (\ref{na}). This and (\ref{sx22}) give 
\[ {\omega _\mathbb{D}}\left( {0,{F_\alpha }} \right) \ge {e^{ - K}}{e^{ - \pi \left( {m\left( {{Q_1}} \right) - 3} \right)}}{e^{ - \pi \left( {2m\left( {{Q^ * }} \right) + 10.52{e^{ - 3\pi }} + 5.26{e^{ - 2\pi m\left( {Q^*} \right)}}} \right)\left( {\sqrt {\log \alpha }  - 1} \right)  - 3\pi \log\alpha }} \]
or
\[\frac{{\log {\omega _\mathbb{D}}{{\left( {0,{F_\alpha }} \right)}^{ - 1}}}}{{\log \alpha }} \le \frac{{K' + \pi \left( {2m\left( {{Q^ * }} \right) + 10.52{e^{ - 3\pi }} + 5.26{e^{ - 2\pi m\left( {Q^*} \right)}}} \right)\left( {\sqrt {\log \alpha }  - 1} \right)  + 3\pi \log \alpha }}{{\log \alpha }},\]
where $K':=K + \pi \left( {m\left( {{Q_1}} \right) - 3} \right)$. Hence taking limits as $\alpha  \to  + \infty$,
\[\mathop {\limsup }\limits_{\alpha  \to  + \infty } \frac{{\log {\omega _\mathbb{D}}{{\left( {0,{F_\alpha }} \right)}^{ - 1}}}}{{\log \alpha }} \le 3\pi. \]
By this and (\ref{liminf}) we take
\[{\tt h}\left( \psi \right) =L= \mathop {\lim }\limits_{\alpha  \to  + \infty } \frac{{\log {\omega _\mathbb{D}}{{\left( {0,{F_\alpha }} \right)}^{ - 1}}}}{{\log \alpha }} = 3\pi. \]
Since $N\left( \alpha  \right) = 1$ for every $\alpha>0$, Corollary \ref{cor} implies that $\mu  = L = 3\pi $.
\\\\
\textbf{Step 4:} Setting
\[{C_1}: = 2\pi \left( {m\left( {{Q^ * }} \right) - 9} \right) > 0,\;{C_2}: = \log \left( {{\pi  \mathord{\left/
			{\vphantom {\pi  8}} \right.
			\kern-\nulldelimiterspace} 8}} \right) + 27\pi  - 6\pi m\left( {{Q^ * }} \right),\]
by (\ref{kla}) we take that for every $\alpha>0$,
\[\frac{{\log {\omega _\mathbb{D}}{{\left( {0,{F_\alpha }} \right)}^{ - 1}}}}{{\log \alpha }} \ge 3\pi  + \frac{{{C_1}}}{{\sqrt {\log \alpha } }} + \frac{{{C_2}}}{{\log \alpha }}.\]
By this and a change of variable, we deduce that
\begin{eqnarray}
\int_1^{ + \infty } {{\alpha ^{3\pi  - 1}}{\omega _\mathbb{D}}\left( {0,{F_\alpha }} \right)d\alpha }  &=& \int_1^{ + \infty } {{\alpha ^{3\pi  - 1 - \frac{{\log {\omega _\mathbb{D}}{{\left( {0,{F_\alpha }} \right)}^{ - 1}}}}{{\log \alpha }}}}d\alpha }  \le \int_1^{ + \infty } {{\alpha ^{ - 1 - \frac{{{C_1}}}{{\sqrt {\log \alpha } }} - \frac{{{C_2}}}{{\log \alpha }}}}d\alpha }   \nonumber \\
&=& \int_1^{ + \infty } {{\alpha ^{ - 1}}{e^{ - \frac{{{C_1}}}{{\sqrt {\log \alpha } }}\log \alpha }}{\alpha ^{{{\log {e^{ - {C_2}}}} \mathord{\left/
 {\vphantom {{\log {e^{ - {C_2}}}} {\log \alpha }}} \right.
 \kern-\nulldelimiterspace} {\log \alpha }}}}d\alpha }  = {e^{ - {C_2}}}\int_1^{ + \infty } {{\alpha ^{ - 1}}{e^{ - {C_1}\sqrt {\log \alpha } }}d\alpha }  \nonumber \\
&=& 2{e^{ - {C_2}}}\int_0^{ + \infty } {{t}{e^{ - {C_1}t}}dt}  = \frac{{2{e^{ - {C_2}}}}}{{{C_1}^2}} <  + \infty. \nonumber
\end{eqnarray}
So, by (\ref{isod}) we infer that $\psi  \in {H^{3\pi }}\left( \mathbb{D} \right)$.
\qed

\end{document}